\documentclass[11pt]{article}
\usepackage{latexsym}
\usepackage{amssymb}
\def\theequation{\thesection.\arabic{equation}}

\font\gorditas = msbm8
\def\bbb#1{\hbox {{\gordas #1}}}
\def\errita{\hbox{\gorditas R}}
                                                                                
\font\gordas = msbm10 at 12pt
\def\bbb#1{\hbox {{\gordas #1}}}
\def\erre{{\bbb R}}

\def\UNO{1\mkern-7mu1}
\def\ee{{\bbb E}}

\def\ce{{\bbb C}}

\def\be{\begin{equation}}

\def\ee{\end{equation}}

\newtheorem{theorem}{Theorem}[section]

\newtheorem{proposition}[theorem]{Proposition}

\newtheorem{remark}[theorem]{Remark}

\begin{document}

\begin{center}
{\bf\Large  
Self-Similar Stable Processes Arising From High-Density\\[.5cm]
 Limits of Occupation Times of Particle Systems}\\[.5cm]
(Self-Similar Stable Processes and Particle Systems)\\[.5cm]
T. BOJDECKI$^{1,*}$, L.~G. GOROSTIZA$^{2,**}$ and A. TALARCZYK$^{3,*}$
\end{center}
\noindent
$^1$ Institute of Mathematics, University of Warsaw, ul. Banacha 2, 02-097
Warszawa,  Poland\\
\phantom{ll} (e.mail: tobojd@mimuw.edu.pl)\\
$^2$ Centro de Investigaci\'on y de Estudios Avanzados, A.P. 14-740, Mexico
07000 D.F.,Mexico\\
\phantom{ll} (e.mail: lgorosti@math.cinvestav.mx)\\
$^3$ Institute of Mathematics, University of Warsaw, ul. Banacha 2, 02-097
Warszawa, Poland\\
\phantom{ll} (e.mail: annatal@mimuw.edu.pl)

\vglue1.5cm
\noindent
{\bf Abstract.}
{\small
We extend results on time-rescaled  occupation time fluctuation limits of the $(d,\alpha, \beta)$-branching particle system $(0<\alpha \leq 2, 0<\beta \leq 1)$ with  Poisson initial condition. The earlier results in the homogeneous case (i.e., with Lebesgue initial intensity measure) were obtained for dimensions $d>\alpha/\beta$ only, since  the particle system becomes locally extinct if $d\le \alpha/\beta$. In this paper we show that by introducing high density of the initial Poisson configuration, limits are obtained for all dimensions, and they coincide with the previous ones if $d>\alpha/\beta$. We also give high-density limits for the systems with finite intensity measures (without high density no limits  exist in this case due to  extinction); the results are different and harder to obtain due to the non-invariance of the measure for the particle motion. In both cases, i.e., Lebesgue and finite intensity measures, for low dimensions ($d<\alpha(1+\beta)/\beta$ and $d<\alpha(2+\beta)/(1+\beta)$, respectively) the limits are determined by non-L\'evy self-similar stable processes. For the corresponding high dimensions the limits are qualitatively different: ${\cal S}'(\erre^d)$-valued L\'evy processes in the Lebesgue case, stable processes constant in time on $(0,\infty)$ in the finite measure case. For high dimensions, the laws of all limit processes are expressed in terms of Riesz potentials. If $\beta=1$, the limits are Gaussian. Limits are also given for particle systems without branching, which yields in particular weighted fractional Brownian motions in low dimensions.
The results are obtained in the setup of weak convergence of ${\cal S}'(\erre^d)$-valued processes.}
\vglue1.5cm
\noindent

\footnote{\kern-.6cm $^*$ Research supported by MNiSW grant 1P03A1129 
(Poland).\\
$^{**}$ Research supported by CONACyT grant 45684-F (Mexico).}
\newpage
\noindent
{\bf 1. Introduction}
\setcounter{section}{1}
\setcounter{equation}{0}
\vglue.5cm

In order to explain the motivations for this paper, we refer briefly to  previous results on occupation times of the  $(d,\alpha,\beta)$-branching particle system, which has been widely studied, and is described as follows.
At time $t=0$ particles are distributed in $\erre^d$ according to a Poisson random measure,  and then they
  evolve  moving and branching independently of each other. The motion is given by the symmetric $\alpha$-stable L\'evy process, $0<\alpha\leq 2$ (called standard $\alpha$-stable process), the lifetime is exponentially distributed with parameter $V$, and the branching law has generating function
\begin{equation}
\label{eq:1.1}
s+\frac{1}{(1+\beta)}(1-s)^{1+\beta},\quad 0<s<1,
\end{equation}
where $0<\beta\leq 1$. This law is critical and belongs to the domain of attraction of a stable law with exponent $1+\beta$. The case $\beta=1$ corresponds to binary branching ($0$ or $2$ particles). 
This is the simplest in a class of branching particle systems that yield essentially the same results.
We also consider the system without branching $(V=0)$.

If the initial particle configuration is given by a homogeneous Poisson random measure, i.e., whose intensity  is the Lebesgue measure $\lambda$, then the system without branching is in equilibrium, the branching system converges towards a non-trivial equilibrium state as time tends to infinity for $d>\alpha/\beta$, and it becomes locally extinct for $d\leq\alpha/\beta$ \cite{13}.

Let $(N_t)_{t\geq 0}$ denote the empirical measure process of the system (with or without branching), i.e., $N_t(A)$ is the number of particles in the set $A\subset \erre^d$ at time $t$. The rescaled occupation time fluctuation process with accelerated time is defined by
\begin{equation}
\label{eq:1.2}
X_T(t)=\frac{1}{F_T}\int^{Tt}_0(N_s-EN_s)ds,\quad t\geq 0,
\end{equation}
where $F_T$ is a suitable norming for convergence as $T\to\infty$. Note that 
 if $\lambda$ is the intensity of the initial Poisson configuration, 
 then $EN_t=\lambda$ for all $t$ due to the invariance of $\lambda$ for the standard $\alpha$-stable process and 
 the criticality of the branching (or no branching). 

With homogeneous Poisson initial condition, functional limit theorems for the process $X_T$ in the branching 
case were obtained in \cite{4,5} for $\beta=1$, where the limit processes are Gaussian, and in \cite{6,7} for $\beta<1$, with $(1+\beta)$-stable limit processes. The limits are dimension-dependent, their main qualitative properties being that for the intermediate dimensions, $\alpha/\beta<d<\alpha(1+\beta)/\beta$, the process has long-range dependence, while for the critical and high dimensions, $d=\alpha(1+\beta)/\beta$ and $d>\alpha(1+\beta)/\beta$, respectively, the processes have independent increments. For high dimensions the limits are ${\cal S}'(\erre^d)$-valued  (${\cal S}'(\erre^d)$ is the space of tempered distributions, the dual of ${\cal S}(\erre^d)$, the space of smooth rapidly decreasing functions), and their laws  are expressed in terms of Riesz potentials. There is  a functional ergodic
 theorem for $d=\alpha/\beta$ \cite{16}. For intermediate dimensions the limit
 has the form $X=K\lambda\xi$, where $K$ is a constant, and $(\xi_t)_{t\geq 0}$ is a real non-L\'evy self-similar $(1+\beta)$-stable process, which for $\beta=1$  is a {\it sub-fractional Brownian motion}, whose properties are described in \cite{3}. 

The first motivation for this paper comes from the fact that 
in the homogeneous case with 
 $\beta=1$ and $d<\alpha$, the covariance of the process $X_T$ has a non-trivial limit as $T\to\infty$, which corresponds to a process $X$ of the same form as above, with a different Gaussian process instead of  sub-fractional Brownian motion, but $X$ is not the limit of $X_T$ because, as recalled above, the particle system becomes locally extinct if $d<\alpha$. Therefore the question arises if it is possible to give a probabilistic meaning (related with the particle system) 
to the process $X$, by taking a different type of limit. Our objective is to show that this can be achieved by letting the density of the initial Poisson configuration tend to infinity in a suitable way as $T\rightarrow \infty$. We will prove a limit theorem for the process $X_T$ for low dimensions, $d<\alpha(1+\beta)/\beta$ (which includes the old intermediate dimensions), and obtain results for the critical and high dimensions as well, by taking an initial Poisson configuration with intensity measure $H_T\lambda$, where $H_T\to\infty$ as $T\to\infty$ (and new normings $F_T$). It turns out that the limits  coincide with the known ones in the cases where the latter exist, i.e., for $d>\alpha/\beta$, and they are new processes for $d\leq\alpha/\beta$, which are also of the form $X=K\lambda\xi$.
For $\beta<1$ and $d<\alpha/\beta,\xi$ is an extension of a  non-L\'evy $(1+\beta)$-stable 
process obtained in \cite{6} for intermediate dimensions (the process in 
\cite{6} has 
 the interesting property that it has two different long-range dependence regimes). For $\beta=1$ and $d<\alpha, \xi$ is a  {\it negative sub-fractional Brownian motion}, which is a real centered Gaussian process with covariance 
\begin{equation}
\label{eq:1.3}
E\xi_s\xi_t=\frac{1}{2}[(s+t)^h+|s-t|^h]-s^h-t^h,\quad s,t\geq 0,
\end{equation}
where $h=3-d/\alpha$.  For $\beta=1$ and $d=\alpha$,  $\xi$ is a  centered Gaussian process with covariance\begin{equation}
\label{eq:1.4}
E\xi_s\xi_t=\frac{1}{2}
\left[
(s+t)^2\log(s+t)
+
(s-t)^2\log|s-t|\right]
-s^2\log s
-t^2\log t
,\quad s,t\geq 0.
\end{equation}
Some properties of these processes are studied in \cite{9}, independently of their origin in particle systems.

Thus, the high-density limits extend the ranges of the parameters of the branching particle system for convergence of $X_T$ obtained  in \cite{4,5,6,7} without high density, so that all cases are now covered, including dimensions $d$ below and at the extinction  border $\alpha/\beta$ .

For completeness, we will also include  high-density limits for the system without branching, but  there are no novelties in the sense that the limits coincide with those for the  homogeneous Poisson case  without high density.

The second  motivation   is the question of what happens with the occupation times of the particle systems if the initial Poisson configuration has finite intensity measure. In this case the branching 
system becomes extinct a.s., while the non-branching system becomes locally extinct a.s. if $d>\alpha$, and if $d\le\alpha$, then 
($1/F_T$)$\int_0^TE\langle N_s,\varphi\rangle ds$ 
 converges to a finite limit for any $\varphi\in{\cal S}(\erre^d)$  and 
($1/F_T$)$\int_0^{Tt}N_sds$ has a non-trivial limit in law (see
 \cite{8}, the latter result   is  akin to the 
Darling-Kac occupation time theorem \cite{11}). For these reasons
 it does not make sense to study asymptotic occupation time  fluctuations. We will show that high density of the initial Poisson condition can be used to compensate  extinction and obtain non-trivial limits for $X_T$. We will consider an initial Poisson configuration with intensity measure $H_T\mu$, where $\mu$ is a finite measure and $H_T\to\infty$ as $T\to\infty$. This yields results for the occupation time fluctuations of the branching and the non-branching  systems, with new types of limits.  These results are different, and significantly more
difficult to obtain  than the previous ones, because the Poisson intensity measure is not invariant for the standard $\alpha$-stable process    
(if the intensity measure is $\mu$, then $EN_t=\mu{\cal T}_t$, where ${\cal T}_t$ is the semigroup of the standard $\alpha$-stable process). 

For the branching system with finite measure $\mu$,  the low, critical and 
high dimensions are $d<\alpha(2+\beta)/(1+\beta)$,  $d=\alpha(2+\beta)/(1+\beta)$, and
 $d>\alpha(2+\beta)/(1+\beta)$, respectively. In the first two cases the limit processes are of the form $K\lambda\xi$. For low dimensions, $\xi$ is a non-L\'evy $(1+\beta)$-stable  process, which is different from the one
 obtained in the homogeneous case. For the critical dimension, $\xi$ is a process  constant in time on $(0,\infty)$, given by a $(1+\beta)$-stable random variable. In these two cases the measure $\mu$ figures only through its total mass, which appears as a constant. For the high dimensions the limit is a process  constant in time on $(0,\infty)$, given by an ${\cal S}'(\erre^d)$-valued $(1+\beta)$-stable random variable whose 
law is expressed by means of a Riesz potential. In this case $\mu$ has a non-trivial effect on the spatial distribution of the limit process. So,  in addition to the  critical borders being different for Lebesgue and finite measures, the limit processes are qualitatively different for the two cases in the corresponding critical and high dimensions.

For the non-branching system, the low, critical and high dimensions for the high-density limits with finite measure are  $d<\alpha, d=\alpha$ and $d>\alpha$, respectively. For $d<\alpha$ the limit has the form $K\lambda\rho$, where $\rho$ is a special case of a {\it weighted fractional Brownian motion} studied in 
\cite{9}, i.e., centered Gaussian 
with covariance
\begin{equation}
\label{eq:1.5}
E\rho_s\rho_z =\int^{s\wedge t}_0 u^{-d/\alpha}[(t-u)^{1-d/\alpha}+(s-u)^{1-d/\alpha}]du,\quad s,t\geq 0.
\end{equation}
For $d=\alpha$ and $d>\alpha$, the limits are constant in time on $(0,\infty)$, analogously to the branching case in the corresponding critical and high dimensions.
They are Gaussian with covariances expressed by means of Riesz potentials.

The  proofs in this paper are analogous to those in \cite{4,5,6,7,8}, but there are new complexities that require a more comprehensive approach.
We will explain the general scheme at the beginning of the proofs, but we stress that its implementation in specific cases is not at all straightforward, and it becomes quite cumbersome technically in the case of finite measure. 
We will refer often to our previous papers (specially \cite{6}) for some technical points, in order to shorten the length of this article, and the main parts of the proofs  given here are devoted  to arguments that involve something new.
The general setting is  weak convergence of ${\cal S}'(\erre^d)$-valued processes, which covers  the cases where the limit process is measure-valued and those where it is ``truly'' ${\cal S}'(\erre^d)$-valued.

We will use the following notions of weak convergence of 
${\cal S}'(\erre^d)$-valued processes (recall that ${\cal S}(\erre^d)$ denotes the space of smooth rapidly decreasing functions, ${\cal S}'(\erre^d)$, the dual of ${\cal S}(\erre^d)$, 
is the space of tempered distributions, and $\langle \cdot, \cdot \rangle$ stands for  duality pairing):

$\Rightarrow_{C}$ is the convergence in law in $C([0,\tau],{\cal S}'(\erre^d))$ for each $\tau>0$;

$\Rightarrow_{C,\varepsilon}$ is the convergence in law in $C([\varepsilon,\tau], {\cal S}'(\erre^d))$ for each $0<\varepsilon<\tau$;

$\Rightarrow_{f}$ is the convergence of finite-dimensional distributions;

$\Rightarrow_{i}$ is the convergence in the integral sense, i.e.,
$X_T\Rightarrow_{i} X$ as $T\to\infty$ if, for any $\tau>0$, the ${\cal S}'(\erre^{d+1})$-random variables $\widetilde{X}_T$ converge in law to $\widetilde{X}$, where $\widetilde{X}$ (and, analogously, $\widetilde{X}_T$) is defined by
\begin{equation}
\label{eq:1.6}
\langle\widetilde{X},\Phi\rangle=\int^\tau_0\langle X(t),\Phi(\cdot,t)\rangle dt,\quad \Phi\in {\cal S}(\erre^{d+1}).
\end{equation}

We denote generic constants by $C, C_1,C_2,\ldots,$ with possible dependencies in parenthesis.
\vglue1cm
\noindent
{\bf 2. Results}
\setcounter{section}{2}
\setcounter{equation}{0}
\vglue.5cm

Before stating the results we introduce two $(1+\beta)$-stable processes
 which appear in the theorems below  $(0<\beta\leq 1$ is fixed).

Let $M$ be the independently scattered $(1+\beta)$-stable measure on $\erre^{d+1}$ with control measure $\lambda_{d+1}$ (Lebesgue measure) and skewness intensity 1, i.e., for each $A\in {\cal B}(\erre^{d+1})$ such that $0<\lambda_{d+1}(A)<\infty, M(A)$ is a $(1+\beta)$-stable random variable with characteristic function
$${\rm exp}\biggl\{-\lambda_{d+1}(A)|z|^{1+\beta}\biggl(1-i({\rm sgn}z)\tan\frac{\pi}{2}(1+\beta)\biggl)\biggr\},\quad z\in\erre,$$
the values of $M$ are independent on disjoint sets, and $M$ is $\sigma$-additive a.s. (see \cite{15}, Definition 3.3.1).

Let $p_t(x)$ denote the transition density of the standard $\alpha$-stable process in $\erre^d$.

We define the following processes:
\begin{eqnarray}
\label{eq:2.1}
\xi_t&=&\int_{\errita^{d+1}}\biggl(\UNO_{[0,t]}(r)\int^t_rp_{u-r}(x)du\biggr)M(drdx),\quad t\geq 0,\\
\label{eq:2.2}
\zeta_t&=&\int_{\errita^{d+1}}\biggl(\UNO_{[0,t]}(r)p_r^{1/(1+\beta)}(x)\int^t_rp_{u-r}(x)du\biggr)M(drdx),\quad t\geq 0,
\end{eqnarray}
where the integral with respect to $M$ is understood in the sense of \cite{15} (3.2-3.4).

\begin{proposition}
\label{p:2.1}
The process $\xi$ is well defined if $d<\alpha(1+\beta)/\beta$, and the process $\zeta$ is well defined if $d<\alpha(2+\beta)/(1+\beta)$.
\end{proposition}

The process $\xi$ is an extension of the one studied in [6].

We denote by ${\cal T}_t$ the semigroup of the standard $\alpha$-stable process, i.e., ${\cal T}_t\varphi=p_t *\varphi$. For $d>\alpha$, we denote by $G$ the potential operator
\begin{equation}
\label{eq:2.3}
G\varphi(x)=\int^\infty_0{\cal T}_t\varphi(x)dt=C_{\alpha,d}\int_{\errita^d}\frac{\varphi(y)}{|x-y|^{d-\alpha}}dy,
\end{equation}
where 
\begin{equation}
\label{eq:2.4}
C_{\alpha, d}=\frac{\Gamma (\frac{d-\alpha}{2})}
{2^\alpha \pi^{d/2}\Gamma (\frac{\alpha}{2})}.
\end{equation}

We start with the high-density branching system described in the Introduction, where the intensity measure of the initial Poisson configuration is  $H_T\lambda$.

\begin{theorem}
\label{t:2.2}
Consider the $(d,\alpha,\beta)$-branching particle system with branching mechanism (\ref{eq:1.1}) and initial intensity $H_T\lambda, H_T\to\infty$. Let $X_T$ be defined by (\ref{eq:1.2}).

\noindent
(a) Assume
\begin{equation}
\label{eq:2.5}
d<\frac{\alpha(1+\beta)}{\beta}.
\end{equation}
Let $H_T$ be such that
\begin{equation}
\label{eq:2.6}
\lim_{T\to\infty}H^{-\beta}_T T^{1-d\beta/\alpha}=0,
\end{equation}
and
\begin{equation}
\label{eq:2.7}
F_T^{1+\beta}=H_TT^{2+\beta-d\beta/\alpha}.
\end{equation}
Then $X_T\Rightarrow_{C}K\lambda\xi$ as $T\to\infty$, where $\xi$ is defined by (\ref{eq:2.1}) and
\begin{equation}
\label{eq:2.8}
K=\biggl(-\frac{V}{1+\beta}\cos\frac{\pi}{2}(1+\beta)\biggr)^{1/(1+\beta)}.
\end{equation}
(b) Assume 
$d=\alpha(1+\beta)/\beta\quad {\it and}\quad F^{1+\beta}_T=H_TT\log T.
$ Then $X_T\Rightarrow_{i}K_1\lambda\eta$ and $X_T\Rightarrow_{f}K_1\lambda\eta$ as $T\to\infty$, where $\eta$ is a real $(1+\beta)$-stable process with stationary independent increments whose distribution is determined by
\begin{equation}
\label{eq:2.9}
E{\rm exp}\{iz\eta_t\}={\rm exp}\biggl\{-t|z|^{1+\beta}\biggl(1-i({\rm sgn} z)\tan \frac{\pi}{2}(1+\beta)\biggr)\biggr\},\quad z\in\erre,\quad t\geq 0,
\end{equation}
and
$$K_1=
\left(- V\int_{\errita^d}\left
(\int^1_0p_r(x)dr\right)^\beta p_1(x)dx\cos \frac{\pi}{2}(1+\beta) \right)^{1/(1+\beta)}.$$
Moreover, if $\beta=1$, the convergence holds in the sense $\Rightarrow_{C}$.

\noindent
(c) Assume 
$d>\alpha(1+\beta)/\beta\quad {\it and}\quad  F^{1+\beta}_T=H_TT.
$ 

\noindent
(i) If $0<\beta<1$, then $X_T\Rightarrow_{i}X$ and $X_T\Rightarrow_{f}X$ as $T\to\infty$, where $X$ is an ${\cal S}'(\erre^d)$-valued $(1+\beta)$-stable process with stationary independent increments whose distribution is determined by
\begin{eqnarray*}
\lefteqn{E{\rm exp}\{i\langle X(t),\varphi\rangle\}}\\
&=&{\rm exp}\biggl\{-K^{1+\beta}t\int_{\errita^d}|G\varphi(x)|^{1+\beta}\biggl(1-i({\rm sgn}G\varphi(x))\tan\frac{\pi}{2}(1+\beta)\biggr)dx\biggr\},\quad \varphi\in{\cal S}(\erre^d),\quad t\geq 0,
\end{eqnarray*}
$K$ is given by (\ref{eq:2.8}) and $G$ by (\ref{eq:2.3}).

\noindent
(ii) If $\beta=1$, then $X_T\Rightarrow_{C}W$ as $T\to\infty$, where $W$ is an 
${\cal S}'(\erre^d)$-valued Wiener process with covariance
$$E(\langle W(s),\varphi_1\rangle\langle W(t),\varphi_2\rangle)=(s\wedge t)
\int_{\errita^d}
\left[V(G\varphi_1(x))(G\varphi_2(x))+2\varphi_1(x)G\varphi_2(x)\right]dx,\quad \varphi_1,\varphi_2\in{\cal S}(\erre^d),$$%
$$\phantom{pppppppppppppppppppppppppppppppppppppppppppppppppppppppppppppppppp}
s,t\geq 0.$$
\end{theorem}

\begin{remark}
\label{eq:r:2.3}

{\rm (a) For $d>\alpha/\beta$, the limits in Theorem 2.2 are exactly the same as in the model without high density  \cite{4,5,6,7}. Thus, if the limits without high density exist, then increasing the initial density of particles does not change the results.

\noindent
(b) Observe  that assumption (\ref{eq:2.6}) is a restriction only if $d<\alpha/\beta$.

\noindent
(c) If $d\leq\alpha/\beta$, then the limit processes are extensions of those studied before \cite{4,5,6,7} in the sense that the ranges of the parameters are increased.

In \cite{6} we discussed some basic properties of $\xi$ defined by (\ref{eq:2.1}) for $\alpha/\beta<d<\alpha(1+\beta)/\beta$. It turns out that $\xi$ has the same properties also for the full ranges of parameters. We collect them in the following proposition.}
\end{remark}
\begin{proposition}
\label{p:2.4}
Assume (\ref{eq:2.5}).

\noindent
(a) $\xi$ is $(1+\beta)$-stable, totally skewed to the right if $\beta<1$.

\noindent
(b) $\xi$ is self-similar with index $b=(2+\beta-d\beta/\alpha)/(1+\beta)$, 
i.e.,
$$(\xi_{at_1},\ldots,\xi_{at_k})\stackrel{d}{=}a^b(\xi_{t_1},\ldots,\xi_{t_k}),\quad a>0.$$

\noindent
(c) $\xi$ has continuous paths.

\noindent
(d) $\xi$ has the long-range dependence property with dependence exponent

\begin{equation}
\label{eq:2.10}
\kappa=\left\{
\begin{array}{lll}
\displaystyle{
\frac{d}{\alpha}}
&{\it if\,\, either}&
\displaystyle{\alpha=2,\,\,{\it or}\,\, \alpha<2\,\,{\it and}\,\,
\beta>\frac{d}{d+\alpha}},\nonumber\\\\
\displaystyle{\frac{d}{\alpha}\biggl(1+\beta-
\frac{d}{d+\alpha}\biggr)}&{\it if}&
\displaystyle{\alpha<2\,\,{\it and}\,\,\beta\leq\frac{d}{d+\alpha}}.
\end{array}\right.
\end{equation}

\end{proposition}
\vglue.5cm
All these properties are obtained  the same way as in \cite{6}. Property (a) follows from the definition, (b) and (c) are consequences of Theorem 2.2, and (d) can be obtained  exactly   as in Theorem 2.7 in \cite{6}. Recall that the {\it dependence exponent} of $\xi$ is defined by
\begin{equation}
\label{eq:2.11}
\kappa =\inf\limits_{z_1z_2\in \errita}\inf\limits_{0\leq u<v<s<t} \sup \{\gamma >0:D_T(z_1, z_2; u,v,s,t)=o(T^{-\gamma})\quad{\rm as}\quad T\rightarrow \infty\},
\end{equation}
where
\begin{eqnarray}
\label{eq:2.12}
\lefteqn{D_T (z_1, z_2; u,v,s,t)}\nonumber\\
 &=& |\log E e^{i(z_1(\xi_v -\xi_u)+z_2(\xi_{T+t}-\xi_{T+s})}
-\log Ee^{iz_1 (\xi_v-\xi_u)}-\log Ee^{iz_2 (\xi_{T+t}-\xi_{T+s})}|,
\end{eqnarray}
see Definition 2.5 in \cite{6}.

The process $\xi$ can be described more explicitly in the case $\beta =1$.

\begin{proposition}
\label{P2.5}
If $\beta =1$ and $d<2\alpha$, then $\xi$ is a centered Gaussian process with covariance
\begin{eqnarray}
 E\xi_s\xi_t\qquad\nonumber\\
&& \kern -2cm =\left\{ 
\begin{array}{lcl}
\displaystyle{\frac{p_1(0)}
{(1-\frac{d}{\alpha})(2-\frac{d}{\alpha})(3-\frac{d}{\alpha})}}
\biggl(\displaystyle{\frac{1}{2}}[(s+t)^{3-d/\alpha}+|s-t|^{3-d/\alpha}]-s^{3-d/\alpha}
-t^{3-d/\alpha}\biggr) &\!\!\!{\it if} &\!\!\! d\neq \alpha,\,\,\,\,\hfil(\theequation)
\label{eq:2.13}\\\\
\addtocounter{equation}{1}
\displaystyle{\frac{p_1(0)}{2}}\biggl(\displaystyle{\frac{1}{2}}[(s+t)^ 2 \log (s+t)+(s-t)^2 \log |s-t|]-s^2 \log s-t^2 \log t\biggr) &\!\!\! {\it if} 
&\!\!\! d=\alpha,\,\,\,\,\hfil(\theequation)
\label{eq:2.14}
\end{array}\right.  \nonumber \addtocounter{equation}{1}\\[.5cm]
 s,t\geq 0. && \nonumber
\end{eqnarray}
\end{proposition}

The Gaussian process $\xi$ with covariance (\ref{eq:2.13}) 
is (up to a multiplicative constant) a sub-fractional Brownian motion if $\alpha <d<2\alpha$, and a negative sub-fractional Brownian motion if $d<\alpha$. These  processes are studied in \cite{3} and \cite{9}, respectively. The latter paper also contains a proof of the non-semimartingale property of the process with covariance (\ref{eq:2.14}).

Next we consider the system without branching. In this case it is known that if the initial intensity measure is $\lambda$, then the limit of $X_T$ exists for all dimensions \cite{4,5}. The observation in Remark 2.3 (a) also applies here, i.e., introducing high density of the initial configuration does not have any effect on the results. For completeness we give the corresponding theorem.

\begin{theorem}
\label{T2.6}
Let $X_T$ be defined by (\ref{eq:1.2}) for a system without branching with initial intensity $H_T\lambda$, $H_T \rightarrow \infty$. 

\noindent
(a) If 
$
d<\alpha \quad{\it and} \quad F_T=H^{1/2}_T T^{1-d/2\alpha}
$,  then $X_T {\Rightarrow_C}K\lambda \vartheta$ as $T\rightarrow \infty$, where $\vartheta$ is a fractional Brownian motion with Hurst parameter $1-d/2\alpha$, i.e., a centered Gaussian process with covariance
$$
E\vartheta_s \vartheta_t=\frac{1}{2}(s^{2-d/\alpha}+t^{2-d/\alpha}-|s-t|^{2-d/\alpha}),\,\, s,t\geq 0,
$$
and
$$
K=\biggl(\frac{2\Gamma (d/\alpha)}
{\pi \alpha (2-d/\alpha)(1-d/\alpha)}\biggr)^{1/2}.
$$
(b) If 
$
d=\alpha \quad {\it and} \quad F_T=(H_TT \log T)^{1/2}
$,  then $X_T{\Rightarrow_C}K_1\lambda \vartheta^{(1)}$ as $T\rightarrow \infty$, where $\vartheta^{(1)}$ is a standard Brownian motion and
$$
K_1=(2^{d-2}\pi^{d/2}d\Gamma (d/2))^{-1/2}.
$$
(c) If 
$
d>\alpha\quad {\it and}\quad  F_T=(H_TT)^{1/2}
$, then $X_T\Rightarrow_CW^{(0)}$ as $T\rightarrow \infty$, where $W^{(0)}$ is an ${\cal S}'(\erre^d)$-valued Wiener valued process with covariance 
$$
E(\langle W^{(0)}(s), \varphi_1\rangle \langle W^{(0)}(t), \varphi_2\rangle) =(s\wedge t)2\int_{\errita^d}\varphi_1 (x) G\varphi_2(x)dx,\,\, s,t\geq 0,
$$
where $G$ is given by (2.3).
\end{theorem}

An analysis of the proofs of Theorems 2.1 in \cite{4} and \cite{5} shows that the same argument can be employed in the present case, therefore we omit the proof of Theorem \ref{T2.6}.

We now pass  to the system with finite initial  intensity measure.
\begin{theorem}
\label{T2.7}
Consider the $(d,\alpha, \beta)$-branching particle system with initial Poisson intensity $H_T\mu$, where $\mu$ is a finite measure and $H_T\rightarrow \infty$. Let $X_T$ be defined by (\ref{eq:1.2}).

\noindent
(a) Assume
\begin{equation}
\label{eq:2.15}
d<\frac{\alpha (2+\beta)}{1+\beta}.
\end{equation}
Let $H_T$ be such that
\begin{equation}
\label{eq:2.16}
\lim\limits_{T\rightarrow \infty} H^{-\beta}_T T=0,
\end{equation}
and
\begin{equation}
\label{eq:2.17}
F^{1+\beta}_T =H_T T^{2+\beta-(d/\alpha)(1+\beta)}.
\end{equation}
Then $X_T\Rightarrow_{C}K\lambda \zeta$ as $T\rightarrow \infty$, where $\zeta$ is defined by (\ref{eq:2.2}) and
$$
K=\biggl(-\frac{V}{1+\beta}\mu (\erre^d)\cos \frac{\pi}{2}(1+\beta)\biggr)^{1/(1+\beta)}.
$$
(b) Assume 
\begin{equation}
\label{eq:2.19}
d={\alpha (2+\beta)\over(1+\beta)},
\end{equation}
let $H_T$ satisfy (\ref{eq:2.16}),
 and 
\begin{equation}
\label{eq:2.20}
F^{1+\beta}_T=H_T \log T.
\end{equation}
Then $X_T\Rightarrow_{C,\varepsilon}K_1 \lambda \eta_1$ as $T\to\infty$, where $\eta_1$ is a $(1+\beta)$-stable random variable, totally skewed to the right (see (\ref{eq:2.9})),
 and
$$
\kern -2cm K_1 =C_{\alpha, d} \biggl(-\frac{V}{1+\beta} \mu (\erre^d) \int_{\errita^d} 
\frac{p_1(y)}
{|y|^{(d-\alpha)(1+\beta)}}dy \cos \frac{\pi}{2}(1+\beta) \biggr)^{1/(1+\beta)},
$$
where $C_{\alpha, d}$ is given by (\ref{eq:2.4}).

\noindent
(c) Assume 
\begin{equation}
\label{eq:2.21}
d>{\alpha (2+\beta)\over(1+\beta)},
\end{equation}
 let $H_T$ satisfy (\ref{eq:2.16})
 and
\begin{equation}
\label{eq:2.22}
F^{1+\beta}_T=H_T.
\end{equation}
(i) If $0<\beta <1$, then $X_T\Rightarrow_{C,\varepsilon}X$ as $T\to\infty$, where $X$ is an ${\cal S}'(\erre^d)$-valued random variable with characteristic function
\begin{equation}
\label{eq:2.23}
Ee^{i\langle X, \varphi\rangle}=\exp \biggl\{-\frac{V}{1+\beta}\int_{\errita^d}|G\varphi (x)|^{1+\beta}\biggl[1-i({\rm sgn}G\varphi (x))\tan \frac{\pi}{2}(1+\beta)\biggr] G\mu (dx) \cos \frac{\pi}{2}(1+\beta)\biggr\},
\end{equation}
where $G$ is given by (2.3).

\noindent
(ii) If $\beta =1$, then $X_T\Rightarrow_{C,\varepsilon}X$ as $T\to\infty$, where $X$ is a centered ${\cal S}'(\erre^d)$-valued Gaussian random variable with covariance
\begin{equation}
\label{eq:2.24}
E(\langle X, \varphi_1 \rangle \langle X, \varphi \rangle) =2 \int_{\errita^d}\biggl[\varphi_1 (x) G\varphi_2 (x)+\frac{V}{2}(G\varphi_1 (x))( G\varphi_2 (x))\biggr]G\mu (dx).
\end{equation}
\end{theorem}
\begin{remark}
\label{R2.8}
\noindent
{\rm (a) In parts (a) and (b) of Theorem \ref{T2.7} the dependence of the limit processes on $\mu$ is quite weak; $\mu (\erre^d)$ appears only in constants. On the other hand, for high dimensions (part (c)) $\mu$ has a non-trivial 
 effect on the  spatial structure of the limit.

\noindent
(b) The limit processes in parts (a) of Theorems \ref{t:2.2} and \ref{T2.7} are similar, while parts (b) and (c) of these theorems (the time structures of the limits) are substantially different. Note also that for $\beta <1$ in the present case the convergence is stronger ($\Rightarrow_{C, \varepsilon}$
 instead of $\Rightarrow_{i}$ and $\Rightarrow_{f}$). On the other hand, it is clear that one cannot expect to have convergence on the whole interval $[0,1]$, since the limit process is discontinuous at $0$.

\noindent
(c) For large dimensions (part (c)), analogously to the case of the Lebesgue measure, the limit for $\beta =1$ is not obtained from (\ref{eq:2.23})
 by putting $\beta =1$. An additional term appears in the covariance, related to the system without branching, due to slower growth of $F_T$ (see (2.26) below).

In the next proposition we collect properties of the process $\zeta$ in  Theorem \ref{T2.7}(a).}
\end{remark}
\begin{proposition}
\label{P2.9}
Assume (\ref{eq:2.15})
 and let $\zeta$ be defined by (\ref{eq:2.2}).

\noindent
(a) $\zeta$ is $(1+\beta)$-stable, totally skewed to the right if $\beta<1$.\\
(b) $\zeta$ is self-similar with index $(2+\beta)/(1+\beta)-d/\alpha$.\\
(c) $\zeta$ has continuous paths.\\
(d) $\zeta$ has long-range dependence exponent $d/\alpha$.

\end{proposition}

The long-range dependence exponent of $\zeta$ does not depend on $\beta$,
whereas the process $\xi$ has two long-range dependence regimes, one depending on $\beta$ (cf. (\ref{eq:2.10})).

\noindent

We remark that the covariance of the Gaussian process $\zeta$ in the case $\beta=1$ does not have a simple form (in contrast with $\xi$, see Proposition \ref{P2.5}).

Finally, we turn to the non-branching high-density system with finite initial intensity measure.

\begin{theorem}
\label{T2.10}
Let $X_T$ be defined by (\ref{eq:1.2}) for a system without branching with initial Poisson intensity $H_T\mu$, where $\mu$ is a finite measure and $H_T\rightarrow \infty$.

\noindent
(a) If $d<\alpha$ and
\begin{equation}
\label{eq:2.25}
F_T =H^{1/2}_T T^{1-d/\alpha},
\end{equation}
then $X_T \Rightarrow_{C} 
(2\mu (\erre^d)/ (1-d/\alpha))^{1/2}p_1(0)\lambda \rho$ as $T\to\infty$,
where $\rho$ is a centered Gaussian process with covariance
\begin{equation}
\label{eq:2.26}
E\rho_s \rho_t=\int^{t\wedge s}_0 u^{-d/\alpha}[(t-u)^{1-d/\alpha}+
(s-u)^{1-d/\alpha}]du,\,\, s,t\geq 0.
\end{equation}
(b) If $d=\alpha$ and $F_T=H^{1/2}_T \log T$, then $X_T\Rightarrow_{C,\varepsilon}(2\mu (\erre^d))^{1/2}p_1 (0)\gamma$ as $T\to\infty$, where $\gamma $ is a standard Gaussian random variable.

\noindent
(c) If $d>\alpha$ and $F_T=H^{1/2}_T$, then $X_T \Rightarrow_{C,\varepsilon}X$ as $T\to\infty$, where $X$ is a centered ${\cal S}'(\erre^d)$-valued Gaussian random variable with covariance
\begin{equation}
\label{eq:2.27}
E(\langle X, \varphi_1 \rangle \langle X, \varphi_2\rangle) =2\int_{\errita^d} \varphi_1 (x) G\varphi_2(x) G\mu (dx),
\end{equation}
with $G$ given by (2.3).
\end{theorem}
\begin{remark}
\label{R2.11} {\rm (a)
 As in the branching case (Theorem \ref{T2.6}), there is a substantial difference in the time structures of the limits for $d\geq \alpha$.

\noindent
(b) The process $\zeta$ with covariance (\ref{eq:2.26})
 belongs to a class of weighted fractional Brownian motions which is discussed in \cite{9}, in particular its long-range dependence is studied.}
\end{remark}
\vglue.5cm
\noindent
{\bf 3. Proofs}
\setcounter{section}{3}
\label{sec:3}
\setcounter{equation}{0}
\vglue.5cm
\noindent
{\bf Proof of Proposition 2.1.} It is known that existence of the  processes $\xi$ and $\zeta$ defined by (\ref{eq:2.1}) and (\ref{eq:2.2}) is equivalent to
\begin{equation}
\label{eq:3.1}
\int_{\errita^d} \int^t_0 \biggl(\int^t_r p_{u-r} (x) dr\biggr)^{1+\beta}drdx <\infty,\quad t\geq 0,
\end{equation}
and
\begin{equation}
\label{eq:3.2}
\int_{\errita^d}\int^t_0 p_r (x) \biggl(\int^t_r p_{u-r}(x) du\biggr)^{1+\beta}drdx<\infty, \quad t\geq 0,
\end{equation}
respectively (see \cite{15}). On the other hand, from Lemma A.1 in \cite{12} it follows that
\begin{equation}
\label{eq:3.3}
\int_{\errita^d} \biggl(\int^t_0 p_u (x)du\biggr)^{1+\beta}dx<\infty, \quad t\geq 0 \quad {\rm if}\quad d<\frac{\alpha (1+\beta)}{\beta},
\end{equation}
and
\begin{equation}
\label{eq:3.4}
\int_{\errita^d} \biggl(\int^t_0 p_u (x) du\biggr)^{2+\beta} dx <\infty, \quad t\geq 0 \quad {\rm if}\quad d<\frac{\alpha (2+\beta)}{1+\beta}.
\end{equation}
(\ref{eq:3.1}) is an immediate consequence of (\ref{eq:3.3}), and (\ref{eq:3.2}) follows from the H\"older inequality and (\ref{eq:3.4}).\hfill $\Box$
\vglue .25cm
\noindent
{\bf General Scheme}
\vglue.5cm

We present a general scheme which will be employed in the convergence proofs. 
We consider a general $(d, \alpha, \beta)$-branching system, initially Poisson
 with intensity measure $\nu_T$.
 Without loss of generality we 
take the time interval $[0,1]$, i.e., $\tau =1$ (see the end of the Introduction). Let $X_T$ be defined by (\ref{eq:1.2}).

Analogously as in \cite{6} (Theorem 2.2) and \cite{7} (Theorem 2.1), we prove that
\begin{equation}
\label{eq:3.5}
\lim_{T\rightarrow \infty} E e^{-\langle \widetilde{X}_T, \Phi\rangle} =
Ee^{-\langle \widetilde{X}, \Phi}\rangle,
\end{equation}
where $X$ is the corresponding limit process, 
$\Phi \in {\cal S}(\erre^{d+1})$, $\Phi \geq 0$, and $\widetilde{X}_T$ and 
$\widetilde{X}$ are defined by (\ref{eq:1.6}). As explained in \cite{6}, due to the special form of the limit (either Gaussian or $(1+\beta)$-stable totally skewed to
the right), (\ref{eq:3.5}) implies $X_T\Rightarrow_i X$. To prove convergence 
$\Rightarrow _C$ (or $\Rightarrow_{C, \varepsilon})$, according to the space-time approach \cite{2}, it suffices to show additionally that the family $\{\langle X_T, \varphi \rangle\}_{T\geq 0}$ is tight in $C([0,1], \erre)$ (or $C([\varepsilon, 1], \erre)$).

For simplicity we consider $\Phi$ of the form
$$
\Phi (x,t)=\varphi\otimes\psi(x,t)
 =\varphi (x) \psi (t), \quad \varphi \in {\cal S}(\erre^d), \psi \in {\cal S}(\erre), \quad \varphi , \psi \geq 0.
$$
It will be clear from the proofs that for general $\Phi$ the argument is analogous.

Denote
\begin{equation}
\label{eq:3.6}
\varphi_T =\frac{1}{F_T} \varphi, \quad \chi (t)=\int^1_t \psi (s) ds, 
\quad \chi_T (t) =\chi \left(\frac{t}{T}\right).
\end{equation}
Let 
\begin{equation}
\label{eq:3.7}
v_T (x,t) =1-E \exp\biggl \{-\int^t_0 \langle N^x_r, \varphi _T\rangle 
\chi_T (T-t+r) dr\biggr\},\,\,0\leq t\leq T,
\end{equation}
where $N^x$ is the empirical process of the branching system started from a single particle at $x$.
It is known that $v_T$ satisfies the equation
\begin{equation}
\label{eq:3.8}
v_T(x,t) =\int^t_0 {\cal T}_{t-u}\biggl[\varphi_T \chi_T (T-u) (1-v_T (\cdot, u))-\frac{V}{1+\beta}v_T^{1+\beta}(\cdot, u)\biggr](x)du, \;0\leq t\leq T,
\end{equation}
(see \cite{6}, (\ref{eq:3.3})). From (\ref{eq:3.7}) and (\ref{eq:3.8}) we obtain immediately

\begin{equation}
\label{eq:3.9}
0\leq v_T \leq 1,
\end{equation}
and
\begin{equation}
\label{eq:3.10}
v_T (x,t)\leq \int^t_0 {\cal T}_{t-u}\varphi_T (x) \chi_T (T-u)du.
\end{equation}
By (\ref{eq:1.2}), the Poisson property, (\ref{eq:3.7}) and $E\langle N^x_t, \varphi\rangle ={\cal T}_t \varphi (x)$, we have

\begin{equation}
\label{eq:3.11}
E e^{-\langle \widetilde{X}_T, \varphi \otimes \psi\rangle} =\exp \biggl\{ -\int_{\errita^d} v_T (x,T) \nu_T (dx) +\int_{\errita^d}\int^T_0 {\cal T}_u \varphi_T (x) \chi_T (u) du \nu_T (dx) \biggr\}.
\end{equation}
Hence, by (\ref{eq:3.8}),
\begin{equation}
\label{eq:3.12}
E e^{-\langle \widetilde{X}_T, \varphi \otimes \psi\rangle}= 
{\rm exp}\left\{{\frac{V}{1+\beta} I_1(T)+I_2(T)-\frac{V}{1+\beta}I_3 (T)}\right\},
\end{equation}
where
\begin{eqnarray}
\label{eq:3.13}
I_1(T) &=& \int_{\errita^d} \int^T_0 {\cal T}_{T-s}
\biggl[\left(\int^s_0 {\cal T}_{s-u}\varphi_T \chi_T (T-u) du\right)^{1+\beta}\biggr] (x) \nu_T (dx),\\
\label{eq:3.14}
I_2 (T) &=& \int_{\errita^d} \int^T_0 {\cal T}_{T-s}
(\varphi_T \chi_T (T-s) v_T (\cdot, s))(x) dx \nu_T (dx),\\
\label{eq:3.15}
I_3 (T) &=& \int_{\errita^d} \int^T_0 {\cal T}_{T-s}\biggl[\left(
\int^s_0 {\cal T}_{s-u} \varphi_T \chi_T (T-u) du\right)^{1+\beta}-
v^{1+\beta}_T (\cdot, s)\biggr](x) ds \nu_T (dx).
\end{eqnarray}
In most of the cases (with the exception of large dimensions and $\beta=1$, where $I_2$ has a nontrivial limit), we prove
\begin{eqnarray}
\label{eq:3.16}
\lim_{T\rightarrow \infty} 
e^{(V/(1+\beta))I_1(T)}&=&
Ee^{-\langle \widetilde{X}, \varphi \otimes \psi\rangle},\\
\label{eq:3.17}
\lim_{T\rightarrow \infty} I_2 (T) &=&0,
\end{eqnarray}
and
\begin{equation}
\label{eq:3.18}
\lim_{T\rightarrow \infty} I_3 (T)=0,
\end{equation}
Note that if $\nu_T =H_T\lambda$,
 then formulas (\ref{eq:3.12})-(\ref{eq:3.14}) have simpler forms due to invariance of $\lambda$ for ${\cal T}_t$. If $\nu_T$
 is finite (hence not invariant under ${\cal T}_t$), then the proofs are more involved.

To prove (\ref{eq:3.17}) we will use the inequality
\begin{equation}
\label{eq:3.19}
I_2(T) \leq \frac{C}{F_T^2} \int_{\errita^d}\int^T_0 \int^T_0 {\cal T}_s (\varphi {\cal T}_u \varphi )(x) duds \nu_T (dx),
\end{equation}
which is an easy consequence of (\ref{eq:3.6}) and (\ref{eq:3.10}).

To obtain (\ref{eq:3.18}) we apply the elementary inequality
$$
(a+b)^{1+\beta}-a^{1+\beta}\leq b^{1+\beta}+(1+\beta) 
a^{{(1+\beta)}/{2}}b^{{(1+\beta)}/{2}},\,\,\, a,b\geq 0, 0<\beta \leq 1,
$$
then by (\ref{eq:3.10}) and (\ref{eq:3.8}) we obtain
\begin{eqnarray*}
0&\leq & I_3 (T) \leq \int_{\errita^d} \int^T_0 {\cal T}_{T-s}
\biggl(\int^s_0 {\cal T}_{s-u} (\varphi _T \chi_T (T-u) v_T)du +
\int^s_0 {\cal T}_{s-u} v^{1+\beta}_T (\cdot, u) du\biggr)^{1+\beta}(x) \nu_T (dx)\\
&+&(1+\beta)\int_{\errita^d} \int^T_0 {\cal T}_{T-s}\biggl[\biggl(\int^s_0 {\cal T}_{s-u} (\varphi_T\chi (T-u)v_T(\cdot, u))+\int^s_0 {\cal T}_{s-u} v^{1+\beta}_T(\cdot, u) du\biggr)^{{(1+\beta)}/{2}}\\
&&\times v_T^{{(1+\beta)}/{2}}(\cdot, s)\biggr](x) \nu_T (dx).
\end{eqnarray*}
We apply the Schwarz inequality to the second term, then we use $(a+b)^{1+\beta}\leq C(a^{1+\beta}+b^{1+\beta})$, $a,b\geq 0$, in both terms, and finally, by (\ref{eq:3.10}), we arrive at
\begin{equation}
\label{eq:3.20}
0\leq I_3 (T)\leq C\biggl(J_1{(T)} +J_2{(T)}+(J_1{(T)}+J{(T)})^{{1}/{2}}
I_1(T)^{{1}/{2}}\biggr),
\end{equation}
where
\begin{eqnarray}
\label{eq:3.21}
J_1(T) &=& \int_{\errita^d} \int^T_0 {\cal T}_{T-s} \biggl[\biggl(
\int^s_0 {\cal T}_{s-u}\left(\varphi_T\int^u_0 {\cal T}_{u-r} \varphi_T dr\right)du\biggr)^{1+\beta}\biggr](x)ds \nu_T (dx)\nonumber\\
&\leq & \frac{1}{F_T^{2+2\beta}}\int_{\errita^d} \int^T_0 {\cal T}_s 
\biggl[\biggr(\int^T_0 {\cal T}_u \left(\varphi \int^T_0 {\cal T}_r 
\varphi dr\biggr)du\right)^{1+\beta}\biggr](x)ds \nu_T (dx),
\end{eqnarray}
and
\begin{eqnarray}
\label{eq:3.22}
J_2 (T) &=& \int_{\errita^d} \int^T_0 {\cal T}_{T-s} 
\biggl[\biggl(\int^s_0 {\cal T}_{s-u} \biggl(
\int^u_0 {\cal T}_{u-r} \varphi_T dr\biggr)^{1+\beta}
du\biggr)^{1+\beta}\biggr](x)ds\nu_T (dx)\nonumber\\
&\leq & \frac{1}{F_T^{(1+\beta) (1+\beta)}}\int_{\errita^d}\int^T_0 {\cal T}_s
\biggl[\biggl(\int^T_0 {\cal T}_u \biggl(\int^T_0 {\cal T}_r \varphi dr\biggr)^{1+\beta}du\biggr)^{1+\beta}\biggr](x)ds \nu_T (dx).
\end{eqnarray}

Given (\ref{eq:3.16}), in order to prove (\ref{eq:3.18}) it suffices to show that 
\begin{equation}
\label{eq:3.23}
\lim_{T\rightarrow \infty}J_1 (T)=0
\end{equation}
and
\begin{equation}
\label{eq:3.24}
\lim_{T\rightarrow \infty} J_2 (T)=0.
\end{equation}

Note that our method of proof of  $\Rightarrow_i$ convergence (based on  equations (\ref{eq:3.8}) and (\ref{eq:3.11})) gives also convergence of finite-dimensional distributions (see, e.g., the proof of Theorem 2.1 in \cite{7}).

In the proofs of tightness of $\{\langle X_{T}, \varphi \rangle\}_{T\geq 2}$ we follow the idea of \cite{6} (proof of Proposition 3.3).
 Fix $0\leq t_1 \leq t_2 \leq 1$ (or $\varepsilon \leq t_1 \leq t_2 \leq 1$ in the proofs of $\Rightarrow_{C,\varepsilon}$ convergence), and let $\psi \in {\cal S}(\erre^d)$ be such that the corresponding $\chi$ (see (\ref{eq:3.6})) satisfies
\begin{equation}
\label{eq:3.25}
0\leq \chi \leq \UNO_{[t_1, t_2]}.
\end{equation} 
We now repeat the argument of the previous part with $\varphi$ replaced by 
$i\theta \varphi$, $\theta >0$. Let $v_{\theta, T}$ be the analogue of 
(\ref{eq:3.7}). Using the inequality
\begin{equation}
\label{eq:3.26}
|1-e^z|\leq 2|z|\quad {\rm if} \quad |e^z |\leq 1, \quad z\in \ce,
\end{equation}
we have
\begin{eqnarray}
\label{eq:3.27}
|v_{\theta, T}(x,t)|&\leq & 2 \theta \int^t_0 \langle N^x_s, \varphi_T\rangle \chi_T (T-t+s)ds\nonumber\\
&=& 2\theta \int^t_0 {\cal T}_{t-s}\varphi_T (x) \chi_T (T-s)ds.
\end{eqnarray}
The function $v_{\theta, T}$ also satisfies equation (\ref{eq:3.8}) with $i\theta\varphi$ 
(we have not assumed $\psi \geq 0$, but it is not needed for (\ref{eq:3.8}) to hold). Hence by (\ref{eq:3.11}) we obtain
\begin{equation}
\label{eq:3.28}
E\exp \{-i\theta \langle \widetilde{X}_T, \varphi \otimes \psi\rangle\}=\exp \{A_\theta (T)+B_\theta (T)\},
\end{equation}
where
\begin{eqnarray}
\label{eq:3.29}
A_\theta(T)&=&i\theta \int_{\errita^d}\int^T_0 {\cal T}_{T-s}
\biggl(\varphi_T \chi_T (T-s) v_{\theta, T}(\cdot, s)\biggr)(x)ds\nu_T (dx),\\
\label{eq:3.30}
B_\theta (T) &=& \frac{V}{1+\beta}\int_{\errita^d} \int^T_0 {\cal T}_{T-s}
\biggl(v^{1+\beta}_{\theta, T}(\cdot, s)\biggr)(x)ds\nu_T (dx).
\end{eqnarray}
From (\ref{eq:3.28}), again by (\ref{eq:3.26}), we have
\begin{equation}
\label{eq:3.31}
0\leq 1-{\rm Re} E \exp \{-i \theta \langle \widetilde{X}_T, \varphi \otimes 
\psi\rangle \}\leq 2(|A_\theta (T)|+|B_\theta (T)|),
\end{equation}
and this implies
\begin{equation}
\label{eq:3.32}
P|\langle \widetilde{X}_T, \varphi \otimes \psi\rangle |\geq \delta)\leq C\delta \int^{{1}/{\delta}}_0 (|A_{\theta}(T)|+|B_\theta (T)|)d\theta, \quad \delta >0,
\end{equation}
(see e.g., \cite{10}, Proposition 8.29). The tightness will be proved if we show that
\begin{equation}
\label{eq:3.33}
|A_\theta (T) |\leq C(\varphi, \sigma, h)\theta^2 (t^h_2-t^h_1)^{1+\sigma},
\end{equation}
and
\begin{equation}
\label{eq:3.34}
|B_\theta (T)|\leq C(\varphi, \sigma, h, V, \beta)\theta^{1+\beta}(t^h_2-t^h_1)^{1+\sigma},
\end{equation}
for some $\sigma, h>0$. Indeed, (\ref{eq:3.32})-(\ref{eq:3.34}) imply, for $0<\sigma <1$,
\be
\label{eq:3.35}
P(|\langle \widetilde{X}_T, \varphi \otimes \psi \rangle|\geq \sigma)\leq 
\frac{C_1}{\delta^2}(t^h_2-t^h_1)^{1+\sigma}.
\ee
We take $\psi$ approximating $\delta_{t_2}-\delta_{t_1}$, and we see that the 
left-hand side of (\ref{eq:3.35}) can be replaced by
$
P(|\langle X_T(t_2),\varphi\rangle -\langle X_T (t_1), \varphi \rangle|\geq \sigma).$
Hence tightness follows by a well-known criterion \cite{1}. (In the case of $\Rightarrow_{C,\varepsilon}$ convergence we use additionally the fact that, as observed above, $\langle X_T (\varepsilon), \varphi \rangle$ converges in law).

In the proofs of (\ref{eq:3.33}) and (\ref{eq:3.34}), we combine (\ref{eq:3.27}) with (\ref{eq:3.29}) or (\ref{eq:3.30}), respectively, obtaining
\begin{eqnarray}
\label{eq:3.36}
|A_\theta (T) | &\leq & 2\theta^2 A(T),\\
\label{eq:3.37}
|B_\theta (T)| &\leq & \frac{2^{1+\beta}V}{1+\beta}\theta^{1+\beta}I_1 (T),
\end{eqnarray}
where
\be
\label{eq:3.38a}
A(T) =\frac{1}{F^2_T}\int_{\errita^d} \int^T_0 \int^s_0 {\cal T}_{T-s}(\varphi {\cal T}_{s-u} \varphi)(x) \chi \left(1-\frac{s}{T}\right)
\chi \left(1-\frac{u}{T}\right)duds\nu_T (dx),
\ee
and $I_1 (T)$ is given by (\ref{eq:3.13}).

Hence we have reduced the proof of tightness to estimating $A(T)$ and $I_1(T)$ by $C(t^h_2-t^h_1)^{1+\sigma}$.

A similar scheme is applied in the cases without branching. We also have (\ref{eq:3.11}) where $v_T$ satisfies (\ref{eq:3.8}) with $V=0$. Then instead of (\ref{eq:3.12}) we have
\be
\label{eq:3.38}
E e^{-\langle \widetilde{X}_T, \varphi \otimes \psi\rangle} =
e^{I\!\!I_1(T)-I\!\!I_2(T)},
\ee
where
\begin{eqnarray}
\label{eq:3.39}
I\!\!I_1 (T) &=& \int_{\errita^d}\int^T_0 \int^s_0 {\cal T}_{T-s} (\varphi_T {\cal T}_{s-u}\varphi_T)(x)\chi_T(T-u)\chi_T(T-s)duds\nu_T (dx),\\
\label{eq:3.40}
I\!\!I_2 (T) &=& \int_{\errita^d} \int^T_0 \int^s_0 {\cal T}_{T-s} (\varphi_T {\cal T}_{s-u} \varphi_T v_T (\cdot, u))(x)\chi_T (T-u)\chi_T (T-s)duds \nu_T (dx),
\end{eqnarray}
and we show that
\be
\label{eq:3.41}
\lim_{T\rightarrow\infty} e^{I\!\!I_1(T)}=Ee^{-\langle \widetilde{X}, \varphi \otimes \psi\rangle},
\ee
and
\be\label{eq:3.42}
\lim_{T\rightarrow \infty}I\!\!I_2 (T)=0.
\ee
Also, the proof of tightness uses the same method as before with 
$B_\theta (T)=0$ (see (\ref{eq:3.30})).

This general scheme is applied in all the proofs (with $\nu_T =H_T\lambda$ or
 $\nu_T=H_T\mu, \mu$ finite measure). However,
 as we have  mentioned in the Introduction,
  its implementation in specific cases is not  
straightforward.
\vglue .25cm
\noindent
{\bf Proof of Theorem 2.2.} We will prove only part (a) of this theorem, as the remaining parts can be obtained the same way as in \cite{5} and \cite{7}. 
Also, since  the proof of (a) is similar to the proof of Theorem 2.2. 
in \cite{6}, we present only the main steps.

We follow the general scheme. Recall that in this case $\nu_T =H_T\lambda$.
 In order to show (\ref{eq:3.16}) it suffices to prove
\be
\label{eq:3.43}
 \lim_{T\rightarrow \infty} I_1 (T) =
\int_{\errita^d} 
\int^1_0 \biggl(\int_{\errita^d}\int^1_s \varphi (y) \psi (r) 
\int^r_s p_{u-s} (x) dudrdy\biggr)^{1+\beta} dsdx,
\ee
and this can be done  the same way as (\ref{eq:3.21}) in \cite{6}. Note that $H_T$ cancels out in $I_1(T)$ (see (\ref{eq:3.13})), and in the proof of 
(\ref{eq:3.21}) in \cite{6},  $\alpha/\beta <d$ was not used, only (\ref{eq:3.1}) 
was  important.

Next, we prove (\ref{eq:3.17}). By (\ref{eq:2.7}),
 after obvious substitutions (\ref{eq:3.19}) has the form
\begin{eqnarray}
\label{eq:3.44}
I_2 (T) &\leq & C H_T^{1-2/(1+\beta)} T^{2(d\beta/\alpha -1)/(1+\beta)}
\int_{\errita^d} \int^1_0 \varphi (x) {\cal T}_{Tu} \varphi (x) dudx\nonumber\\
&\leq & C_1 T^{2(d\beta/\alpha-1)/(1+\beta)-1} \int_{\errita^d} 
\frac{1-e^{-T|x|^\alpha}}{|x|^\alpha}|\widehat{\varphi} (x)|^2 dx,
\end{eqnarray}
where we have used $1-2/(1+\beta)\leq 0$, the Plancherel formula, and the fact that $\widehat{{\cal T}_u\varphi}(x)=e^{-u|x|^\alpha}\widehat{\varphi}(x)$ ($\widehat{\,\,}$ denotes Fourier transform, defined by $\widehat{\varphi}(z)=\int_{\errita^d}e^{ix\cdot z}\varphi(x)dx, z\in \erre^d$, where $\cdot$ is the scalar product in $\erre^d$). 
Hence it is clear that (\ref{eq:3.17}) holds if 
$\alpha /\beta <d<\alpha 
(1+\beta)/\beta$ and if $d<\alpha/\beta$ (we use  
$(1-e^{-T|x|^\alpha})/(T|x|^\alpha)\leq C$).

For $d=\alpha/\beta$, we estimate the right-hand side of (\ref{eq:3.44})
 by
$$
C_1 T^{-{1}/{2}}\int_{\errita^d}\biggl(\frac{1-e^{-T|x|^\alpha}}{T|x|^\alpha}\biggr)^{{1}/{2}}\frac{1}{|x|^{\alpha/2}}|\widehat{\varphi}(x)|^2 dx,
$$
which tends to $0$ as $T\rightarrow \infty$, since $\alpha \leq d$.

To prove (\ref{eq:3.18}) we show (\ref{eq:3.23}) and (\ref{eq:3.24}). By (\ref{eq:2.7}), on the right-hand side of (\ref{eq:3.21}) $H_T$ appears only as a factor $H_T/H^2_T$ 
(which is bounded), and the remaining term tends to $0$ by the same argument as in \cite{6} (see the proof of (\ref{eq:3.33}) therein, where only (\ref{eq:3.1}) was used). Hence we obtain (\ref{eq:3.23}).

So far we have not used the assumption (\ref{eq:2.6});
 it will be needed in the proof of (\ref{eq:3.24}).

By (\ref{eq:3.22}), repeating the argument of \cite{6} (see (\ref{eq:3.35}) therein and the estimates following it), we obtain
$$
J_2 (T) \leq CH^{-\beta}_T T^{1-{d}{\beta}/\alpha}\rightarrow 0,
$$
by assumption (\ref{eq:2.6}). This completes the proof of (\ref{eq:3.5}) by (\ref{eq:3.43})
 and (\ref{eq:2.1}).

In order to prove tightness, we show (\ref{eq:3.33}) and 
(\ref{eq:3.34}) with $h=1$ and
\begin{equation}
\label{eq:3.45}
0<\sigma <\biggl(1+\beta -\frac{d\beta}{\alpha}\biggr)\wedge \beta.
\end{equation}
Note that in (\ref{eq:3.13}) $H_T$ cancels out, and then the proof of (\ref{eq:3.34}) follows the lines of the proof of (\ref{eq:3.48}) in \cite{6}. The assumption $\sigma <\beta$ is needed in order to have $(1+\beta)/(1+\sigma) >1$ (see (\ref{eq:3.55}) in \cite{6}).

It remains to show (\ref{eq:3.33}). By (\ref{eq:3.38a}) we have
\begin{eqnarray}
A(T)&=&\frac{1}{(2\pi)^d}
\frac{H_T}{H_T^{2/(1+\beta)}}
\frac{T^2}{T^{2(2+\beta-d\beta/\alpha)/(1+\beta)}}\int^1_0\chi(s)\int_{\errita^d}|\widehat{\varphi}(x)|^2\int^1_se^{-T(u-s)|x|^\alpha}\chi(u)dudxds\nonumber\\
\label{eq:3.46}
&\leq&\frac{1}{(2\pi)^d}T^{-2(1-d\beta/\alpha)/(1+\beta)}(t_2-t_1)^{1+\sigma}\int_{\errita^d}|\widehat{\varphi}(x)|^2\left(\frac{1}{1+T|x|^\alpha}\right)^{1-\sigma}dx,
\end{eqnarray}
where in the last estimate we used
\begin{eqnarray}
\int^1_se^{-T(u-s)|x|^\alpha}\chi(u)du&\leq&\biggl(\int^1_se^{-T(u-s)|x|^\alpha}dr\biggr)^{1-\sigma}\biggl(\int^1_s\chi(u)du\biggr)^\sigma\nonumber\\
\label{eq:3.47}
&\leq&\biggl(\frac{1}{1+T|x|^\alpha}\biggr)^{1-\sigma}(t_2-t_1)^\sigma
\end{eqnarray}
for any $0<\sigma\leq 1$. Hence (\ref{eq:3.33}) follows immediately if $d\leq \alpha/\beta$. For $d>\alpha/\beta$ (this case was also treated in \cite{6}) 
we write $(1+T|x|^\alpha)^{-(1-\sigma)}\leq T^{\sigma-1}
|x|^{\alpha(\sigma-1)}$, we use $\alpha<d$ and (\ref{eq:2.5}),
 and we see that for $\sigma$ satisfying (\ref{eq:3.45}) the estimate (\ref{eq:3.33}) holds since the term involving $T$ tends to $0$. \hfill$\Box$
\vglue.5cm
\noindent
{\bf Proof of Proposition 2.5} From (\ref{eq:2.1}), for $\beta=1$ we have
\begin{eqnarray}
E\xi_s\xi_t&=&\int^{s\wedge t}_0\int_{\errita^d}\int^s_r\int^t_rp_{u-r}(x)p_{u'-r}(x)du'dudxdr\nonumber\\
\label{eq:3.48}
&=&p_1(0)\int^{s\wedge t}_0\int^s_r\int^t_r(u+u'-2r)^{-d/\alpha}du'dudr,
\end{eqnarray}
by the Chapman-Kolmogorov identity and the self-similarity of the standard $\alpha$-stable process. Hence (\ref{eq:2.13})
 and (\ref{eq:2.15}) follow by calculus. \hfill$\Box$
\vglue.5cm
\noindent
{\bf Proof of Theorem 2.7.}

{\it Proof of part (a)}. According to the general scheme, we show (\ref{eq:3.16}), which amounts to proving
\begin{equation}
\label{eq:3.49}
\lim_{T\to\infty}I_1(T)=\mu(\erre^d)\int^1_0\int_{\errita^d}p_s(y)\biggl(\int^1_sp_{u-s}(y)\chi(u)du\biggr)^{1+\beta}dyds\biggl(\int_{\errita^d}\varphi(z)dz\biggr)^{1+\beta}.
\end{equation}

In (\ref{eq:3.13}) with $\nu_T=H_T\mu$
 we substitute $u'=(T-u)/T, s'=(T-s)/T$, we use the self-similarity  of the $\alpha$-stable density,
\begin{equation}
\label{eq:3.50}
p_{st}(x)=t^{-d/\alpha}p_s(xt^{-1/\alpha}),
\end{equation}
and by (\ref{eq:2.17}) and (\ref{eq:3.6}) we obtain
\begin{eqnarray}
\kern-1cm I_1(T)&=&T^{-d/\alpha}\int_{\errita^d}\int^1_0\int_{\errita^d}p_s((x-y)T^{-1/\alpha})\nonumber\\
&&\times\left(\int^1_s\int_{\errita^d}p_{u-s}((y-z)T^{-1/\alpha})\varphi(z)\chi(u)dzdu\right)^{1+\beta}dyds\mu(dx)\nonumber\\
\label{eq:3.51}
&=&\int_{\errita^d}\int^1_0\int_{\errita^d}p_s(xT^{-1/\alpha}-y)\biggl(\int^1_s\int_{\errita^d}p_{u-s}(y-z)\widetilde{\varphi}_T(z)\chi(u)dzdu\biggr)^{1+\beta}dyds\mu(dx),
\end{eqnarray}
where
\begin{equation}
\label{eq:3.52}
\widetilde{\varphi}_T(z)=T^{d/\alpha}\varphi(zT^{1/\alpha}).
\end{equation}
We denote
\begin{equation}
\label{eq:3.53}
h_s(y)=\int^1_sp_{u-s}(y)\chi(u)du,
\end{equation}
and we write
\begin{equation}
\label{eq:3.54}
I_1(T)=I'_1(T)+I''_1(T),
\end{equation}
where
\begin{eqnarray}
\label{eq:3.55}
I'_1(T)&=&\int_{\errita^d}\int^1_0p_s*h_s^{1+\beta}(xT^{-1/\alpha})ds\mu(dx)
\biggl(\int_{\errita^d}\varphi(z)dz\biggr)^{1+\beta},
\phantom{dddddddddeeeeeddddddddeeeeeeee}\\
\label{eq:3.56}
I''_1(T)&=&\int_{\errita^d}\int^1_0\int_{\errita^d}p_s(xT^{-1/\alpha}-y)\biggl[(h_s*\widetilde{\varphi}_T(y))^{1+\beta}-
\left(h_s(y)\int_{\errita^d}\varphi(z)dz\right)^{1+\beta}\biggr]dydx\mu(dx).
\end{eqnarray}
By (\ref{eq:3.4}), it is not difficult to see that $I'_1(T)$ converges to the right-hand side of 
(\ref{eq:3.49}) . Therefore, to obtain (\ref{eq:3.49}) it suffices to show that $\lim_{T\to\infty}I''_1(T)=0$.
Fix any $\delta$ satisfying
\begin{equation}
\label{eq:3.57}
\frac{d}{\alpha}-\frac{1}{1+\beta}<\delta<1,
\end{equation}
(such $\delta$ exists by (\ref{eq:2.15})).
 We estimate (\ref{eq:3.56}) applying the H\"older inequality to the integrals with respect to the measure $dys^{-\delta}ds\mu(dx)$, obtaining
\begin{eqnarray}
\lefteqn{|I''_1(T)|}\nonumber\\
&\leq&\biggl(\int_{\errita^d}
\int^1_0s^{-\delta}\int_{\errita^d}
\left(s^\delta p_s(xT^{-1/\alpha}-y)\right)^{2+\beta}dyds
\mu(dx)\biggr)^{1/(2+\beta)}\nonumber\\
\label{eq:3.58}
&&\times\biggl(\int_{\errita^d}\int^1_0s^{-\delta}\int_{\errita^d}\left|(h_s*\widetilde{\varphi}_T(y))^{1+\beta}-
\left(h_s(y)\int_{\errita^d}\varphi(z)dz\right)^{1+\beta}\right|^{(2+\beta)/(1+\beta)}dyds\mu(dx)\biggr)^{(1+\beta)/(2+\beta)}.\nonumber\\
&&
\end{eqnarray}
The first factor does not depend on $T$ and is finite by (\ref{eq:3.50}),
(\ref{eq:3.57}) and  finiteness of $\mu$.

By (\ref{eq:3.4}) and the form of $\widetilde{\varphi}_T$ (see (\ref{eq:3.52})) $(h_s*\widetilde{\varphi}_T)^{1+\beta}$ converges to $(h_s\int_{\errita}\varphi(z)dz)^{1+\beta}$ in \\
$L^{(2+\beta)/(1+\beta)}(\erre^d)$ for any $s\in [0,1]$. Moreover, $h_s(y)\leq \int^1_0p_u(y)du$
 (see 
(\ref{eq:3.53})), hence it is not hard to see that the dominated convergence theorem can be applied to show that the right-hand side of (\ref{eq:3.58}) tends to $0$ as $T\to\infty$. So (\ref{eq:3.49})  is proved, and therefore so is
 (\ref{eq:3.16}).

To show (\ref{eq:3.17}) we make obvious substitutions in the right-hand side of (\ref{eq:3.19}) and use self-similarity, obtaining
\begin{equation}
\label{eq:3.59}
I_2(T)\leq C\frac{H_TT^{2-2d/\alpha}}{F^2_T}
\int_{\errita^d}\int_{\errita^{2d}}f(xT^{1/\alpha}-y)f(y-z)\widetilde{\varphi}_T(y)\widetilde{\varphi}_T(z)dzdy\mu(dx),
\end{equation}
where $\widetilde{\varphi}_T$ is given by (\ref{eq:3.52}), and
\begin{equation}
\label{eq:3.60}
f(x)=\int^1_0p_s(x)ds.
\end{equation}

By the H\"older inequality applied to the integral on $z,y$, we have
$$I_2(T)\leq C\frac{H_TT^{2-2d/\alpha}}{F^2_T}\mu(\erre^d)||f||^2_{2+\beta}||\widetilde{\varphi}_T||^2_{(2+\beta)/(1+\beta)}.$$
(\ref{eq:3.4}),(\ref{eq:3.52}) 
and (\ref{eq:2.17}) imply
$$I_2(T)\leq C_1T^{2(d/\alpha(2+\beta)-1/(1+\beta))}\to 0,$$
by (\ref{eq:2.15}).

To complete the proof of (\ref{eq:3.5}) we show (\ref{eq:3.23}) and (\ref{eq:3.24}). From (\ref{eq:3.21}), by a similar argument as in (\ref{eq:3.59}) we obtain
$$J_1(T)\leq\frac{H_TT^{1+2(1+\beta)-2(d/\alpha)(1+\beta)}}{F_T^{2(1+\beta)}}\int_{\errita^d}f*(f*\widetilde{\varphi}_T(f*
\widetilde{\varphi}_T))^{1+\beta}(xT^{-1/\alpha})\mu(dx),$$
with $f,\widetilde{\varphi}_T$ as above. Applying the H\"older and Young inequalities several times we obtain
$$||f*(f*\widetilde{\varphi}_T(f*\widetilde{\varphi}_T))^{1+\beta}||_\infty\leq ||f||^{3+\beta}_{2+\beta}||\widetilde{\varphi}_T||^{1+\beta}_1||\widetilde{\varphi}_T||^{1+\beta}_{(2+\beta)/(1+\beta)}.$$
Hence, by (\ref{eq:2.17}), (\ref{eq:3.52}), (\ref{eq:3.60}) and (\ref{eq:3.4}),
$$J_1(T)\leq C T^{(d/\alpha)(1+\beta)/(2+\beta)-1}\to 0,$$
by (\ref{eq:2.15}).

Finally, by (\ref{eq:3.22}) and the usual argument we get
$$J_2(T)\leq \frac{H_TT^{2+\beta+(1+\beta)(1+\beta)-(d/\alpha)(1+\beta)(1+\beta)}}{F_T^{(1+\beta)(1+\beta)}}
\int_{\errita^d}f*(f*(f*\widetilde{\varphi}_T)^{1+\beta})^{1+\beta}
(xT^{-1/\alpha})\mu(dx).$$
In this case
$$||f*(f*(f*\widetilde{\varphi}_T)^{1+\beta})^{1+\beta}||_\infty\leq ||f||^{2+\beta}_{2+\beta}||f||^{(1+\beta)(1+\beta)}_{1+\beta}
||\widetilde{\varphi}_T||^{(1+\beta)(1+\beta)}_1\leq C,$$
since $||\widetilde{\varphi}_T||_1=||\varphi||_1$. Hence, by (\ref{eq:2.17})and (2.16),
$$J_2(T)\leq C\frac{T}{H_T^\beta}\to 0.$$

We now pass to the proof of tightness.
To prove (\ref{eq:3.33}) we rewrite (\ref{eq:3.38a}) as
$$A(T)=
\frac{H_TT^2}{F^2_T}
\int^1_0\int^1_s\int_{\errita^d}\varphi(x){\cal T}_{T(u-s)}
\varphi(x)(\mu{\cal T}_{Ts})(dx)\chi(u)\chi(s)duds.$$
We use the following identity, which holds for any finite measure $m$,
$$\int_{\errita^d}\varphi_1(x)\varphi_2(x)m(dx)=
\frac{1}{(2\pi)^{2d}}\int_{\errita^{2d}}\widehat{\varphi}_1(x)\widehat{\varphi}_2(y)
\overline{\widehat{m}(x+y)}dxdy,$$
obtaining
\begin{equation}
\label{eq:3.61} 
A(T)=\frac{H_TT^2}{(2\pi)^{2d}F^2_T}\int^1_0\int^1_s\int_{\errita^{2d}}
\widehat{\varphi}(x)e^{-T(u-s)|y|^\alpha}\widehat{\varphi}(y)
e^{-Ts|x+y|^\alpha}\overline{\widehat{\mu}(x+y)} dxdy\chi(u)\chi(s)duds.
\end{equation}
Fix $h$ satisfying
\begin{equation}
\label{eq:3.62}
\left(1-\frac{d}{\alpha}\right)^+<h<\biggl(\frac{2+\beta}{1+\beta}-\frac{d}{\alpha}\biggr)\wedge 1.
\end{equation}
The function $r\to r^{1-h}e^{-r}$ is bounded on $[0,\infty)$, hence we have
from (\ref{eq:3.61})
\begin{eqnarray*}
\lefteqn{A(T)}\\
&\leq&C\frac{H_TT^2}{F^2_TT^{2(1-h)}}
\int^1_0\int^1_s(u-s)^{h-1}s^{h-1}\chi(u)\chi(s)duds
\int_{\errita^{2d}}|\widehat{\varphi}(x)\widehat{\varphi}(y)||y|^{\alpha(h-1)}|x+y|^{\alpha(h-1)}dxdy\\
&\leq&C_1\frac{T^{-2(2+\beta)/(1+\beta)+2d/\alpha+2h}}
{H_T^{(1-\beta)/(1+\beta)}}\int^{t_2}_{t_1}\int^{t_2}_s(u-s)^{h-1}s^{h-1}duds,
\end{eqnarray*}
by (\ref{eq:2.17}) , (\ref{eq:3.25}), and since $\alpha(1-h)<d$ by 
(\ref{eq:3.62}). The
right-hand side of (\ref{eq:3.62}) implies that the  term involving $T$ is bounded, so it is easy to see that (\ref{eq:3.33})  is obtained with $\sigma=h$.

In order to prove (\ref{eq:3.34}) we use (\ref{eq:3.37}). By (\ref{eq:3.51}) 
and (\ref{eq:3.25})  we have
$$I_1(T)\leq \int_{\errita^d}
\left[R_1(xT^{-1/\alpha})+R_2(xT^{-1/\alpha})\right]\mu(dx),$$
where
\begin{eqnarray}
\label{eq:3.63}
R_1(x)&=&\int^{t_1}_0\int_{\errita^d}p_s(x-y)
\left(\int^{t_2}_{t_1}\int_{\errita^d}p_{u-s}(y-z)
\widetilde{\varphi}_T(z)dzdy\right)^{1+\beta}dyds,\\
\label{eq:3.64}
R_2(x)&=&\int^{t_2}_{t_1}\int_{\errita^d}p_s(x-y)\left(\int^{t_2}_s\int_{\errita^d}p_{u-s}(y-z)\widetilde{\varphi}_T(z)dzdy\right)^{1+\beta}dyds.
\end{eqnarray}
Since $\mu$ is finite, it is enough to show that
\begin{equation}
\label{eq:3.65}
\sup_{x\in\errita^d}R_j(x)\leq C(t^h_2-t^h_1)^{1+\sigma},\,\, j=1,2
\end{equation}
for some positive $h$ and $\sigma$.

Fix $\delta>0$ satisfying (\ref{eq:3.57}) and
\begin{equation}
\label{eq:3.66}
\delta\beta >\frac{(1+\beta)^2}{2+\beta}\frac{d}{\alpha}-1.
\end{equation}
(\ref{eq:3.66})  holds for $\delta$ sufficiently close to 1 because from 
(\ref{eq:2.15})  it follows that
$$\frac{(1+\beta)^2}{2+\beta}\frac{d}{\alpha}-1<\beta.
$$
For any fixed $s\in[0,t_1]$, by the Jensen inequality applied to the measure
$$\frac{(u-s)^{-\delta}}{\int^{t_2}_{t_1}(r-s)^{-\delta}dr}\UNO_{[t_1,t_2]}(u)du$$
(this
 trick is borrowed from \cite{12}), we have
\begin{eqnarray}
R_1(x)&\leq&\kern-.3cm \int^{t_1}_0\int_{\errita^d}p_s(x-y)\biggl(
\int^{t_2}_{t_1}(r-s)^{-\delta}dr\biggr)^\beta\nonumber\\
&&\times\int^{t_2}_{t_1}(u-s)^{-\delta}\left(\int_{\errita^d}(u-s)^\delta p_{u-s}(y-z)\widetilde{\varphi}_T(z)dz\right)^{1+\beta}dudyds\nonumber\\
\label{eq:3.67}
&\leq&\kern-.3cm C(t_2-t_1)^{(1-\delta)\beta}\int^{t_2}_{t_1}\int^{t_1}_0(u-s)^{\delta\beta}\int_{\errita^d}p_s(x-y)(p_{u-s}*\widetilde{\varphi}_T(y))^{1+\beta}dydsdu.
\end{eqnarray}
By the H\"older and Young inequalities,
\begin{eqnarray}
\int_{\errita^d}\ldots dy&\leq&\left(
\int_{\errita^d}p^{2+\beta}_s(y)dy\right)^{1/(2+\beta)}\left(\int_{\errita^d}p_{u-s}^{2+\beta}(y)dy\right)^{(1+\beta)/(2+\beta)}\left(\int_{\errita^d}\widetilde{\varphi}_T(z)dz\right)^{1+\beta}\nonumber\\
\label{eq:3.68}
&=&Cs^{-(d/\alpha)(1+\beta)/(2+\beta)}(u-s)^{-(d/\alpha)(1+\beta)^2/(2+\beta)},
\end{eqnarray}
where we have used (\ref{eq:3.50})  and (\ref{eq:3.52}). Observe that by 
(\ref{eq:2.15}) and (\ref{eq:3.66}), 
$$1-\frac{d}{\alpha}\frac{1+\beta}{2+\beta}>0 \quad{\rm and}\quad
1+\delta\beta-\frac{d}{\alpha}\frac{(1+\beta)^2}{2+\beta}>0.$$ 
Hence, combining 
(\ref{eq:3.68})  with (\ref{eq:3.67}), substituting $s'=s/u$ and estimating the integral on $s'$
  by the corresponding value of the beta function, 
\begin{eqnarray}
R_1(x)&\leq&C_1(t_2-t_1)^{(1-\delta)\beta}\int^{t_2}_{t_1}
u^{1-(d/\alpha)(1+\beta)/(2+\beta)+\delta\beta-(d/\alpha)(1+
\beta)^2/(2+\beta)}du\nonumber\\
\label{eq:3.69}
&\leq&C_2(t^{h'}_2-t^{h'}_1)^{1+(1-\delta)\beta},
\end{eqnarray}
where
$$h'=\min\{1,2{-(d/\alpha)(1+\beta)/(2+\beta)+\delta\beta-(d/\alpha)
(1+\beta)^2/(2+\beta)}\}.$$

To estimate $R_2$ (see (\ref{eq:3.64})) we use the H\"older inequality as in 
(\ref{eq:3.58}), and then the Young inequality, obtaining
\begin{eqnarray}
R_2(x)&\leq&\left[\int^{t_2}_{t_1}s^{-\delta}\int_{\errita^d}(s^\delta p_s(x-y))^{2+\beta}dyds\right]^{1/(2+\beta)}\nonumber\\
&&\times\left[\int^{t_2}_{t_1}s^{-\delta}\int_{\errita^d}\left(\int^{t_2-t_1}_0\int_{\errita^d}p_u(y-z)\widetilde{\varphi}_T(z)dzdu\right)^{2+\beta}dyds\right]^{(1+\beta)/(2+\beta)}\nonumber\\
&=&C\biggl(\int^{t_2}_{t_1}s^{(\delta-d/\alpha)(1+\beta)}ds
\biggr)^{1/(2+\beta)}
(t^{1-\delta}_2-t^{1-\delta}_1)^{(1+\beta)/(2+\beta)}\nonumber\\
&&\times\biggl[\int_{\errita^d}\left(\int^{t_2-t_1}_0p_u*\widetilde{\varphi}_T(y)du\right)^{2+\beta}dy\biggr]^{(1+\beta)/(2+\beta)}\nonumber\\
\label{eq:3.70}
&\leq&C_1(t^{h''}_2-t^{t^{h''}}_1)Q^{(1+\beta)/(2+\beta)},
\end{eqnarray}
where
$$h''=\min\left\{1-\delta,1+\left(\delta-\frac{d}{\alpha}
\right)(1+\beta)\right\}$$
(note that $h''>0$ by (\ref{eq:3.57})), and
$$Q=\int_{\errita^d}\left(\int^{t_2-t_1}_0p_u(y)du\right)^{2+\beta}dy.$$
To estimate $Q$ we substitute $u'=u/(t_2-t_1)$, we use self-similarity and  
(\ref{eq:3.4}), obtaining
$$Q=C(t_2-t_1)^{2+\beta-(d/\alpha)(1+\beta)},$$
the exponent  being positive by  (\ref{eq:2.15}).
Combining this with (\ref{eq:3.70})  we have
$$R_2(x)\leq C_2(t^{h''}_2-t^{h''}_1)^{2+\beta-(d/\alpha)(1+\beta)^2/(2+\beta)}.$$
This and (\ref{eq:3.69})  imply (\ref{eq:3.65})  with
$$h=\min\{h',h''\}\quad{\rm and}\quad \sigma=\min 
\left\{(1-\delta)\beta,1+\beta-\frac{d}{\alpha}\frac{(1+\beta)^2}{2+\beta}\right\}.$$

This proves (\ref{eq:3.34})  and completes the proof of part (a) of the theorem.

{\it Proof of part (b).} According to the general scheme,
 we prove 
(\ref{eq:3.16}), and it is easy to see that to this end it suffices to show that
\begin{equation}
\label{eq:3.71}
\lim_{T\to\infty}I_1(T)=\mu(\erre^d)C^{1+\beta}_{\alpha,d}
\int_{\errita^d}p_1(y)|y|^{-(d-\alpha)(1+\beta)}dy
\left(\int_{\errita^d}\varphi(z)dz\right)^{1+\beta}(\chi(0))^{1+\beta}.
\end{equation}
Observe that (\ref{eq:2.19})  implies that $d>\alpha$, hence
\begin{equation}
\label{eq:3.72}
\sup_{x\in\errita^d}G\varphi(x)<\infty,
\end{equation}
where $G$ is defined by (\ref{eq:2.3}). This fact implies in particular that
\begin{equation}
\label{eq:3.73}
\lim_{T\to\infty}I_1(T)=\lim_{T\to\infty}I'_1(T),
\end{equation}
where
$$I'_1(T)=\frac{1}{\log T}\int_{\errita^d}\int^{T-1}_0{\cal T}_{T-s}\left(
\int^s_0{\cal T}_{s-u}\varphi \chi\left(\frac{T-u}{T}\right)du\right)^{1+\beta}(x)ds\mu(dx)$$
(see (\ref{eq:3.13}), (\ref{eq:3.6}) and (\ref{eq:2.20})). By obvious substitutions,
\begin{eqnarray*}
I'_1(T)&=&\frac{1}{\log T}\int_{\errita^d}\int^T_1\int_{\errita^d}p_s(x-y)\left(\int^{T-s}_0\int_{\errita^d}p_u(y-z)\varphi(z)\chi\left(\frac{u}{T}+\frac{s}{T}\right)dzdu\right)^{1+\beta}dyds\mu(dx)\\
&=&\frac{1}{\log T}\int_{\errita^d}\int^T_1\int_{\errita^d}p_1(xs^{-1/
\alpha}-y)\\
&&\times\left(\int^{T-s}_0\int_{\errita^d}s^{-d/\alpha}p_{u/s}(y-zs^{-1/\alpha})\varphi(z)\chi\left(\frac{u}{T}+\frac{s}{T}\right)dzdu\right)^{1+\beta}dyds\mu(dx),
\end{eqnarray*}
where we have used self-similarity and the substitution $y'=ys^{-1/\alpha}$. Next, we substitute $u'=u/s$, and using (\ref{eq:2.19})  we get
\begin{eqnarray*}
I'_1(T)&=&\frac{1}{\log T}\int_{\errita^d}\int^T_1\int_{\errita^d}s^{-1}p_1(xs^{-1/\alpha}-y)\\
&&\times\left(\int^{T/s-1}_0\int_{\errita^d}p_u(y-zs^{-1/\alpha})\varphi(z)
\chi\left(\frac{us}{T}+\frac{s}{T}\right)drdu\right)^{1+\beta}dyds\mu(dx).
\end{eqnarray*}
Now we make the substitution $s'=\log s/\log T$, which is the main trick in calculating the limit.
We obtain
\begin{eqnarray}
I'_1(T)&=&\int_{\errita^d}\int^1_0\int_{\errita^d}p_1(xT^{-s/\alpha}-y)\nonumber\\
\label{eq:3.74}
&&\times\left(\int^{T^{1-s}-1}_0\int_{\errita^d}p_u(y-zT^{-s/\alpha})\varphi(z)\chi((u+1)T^{s-1})dzdu\right)^{1+\beta}dyds\mu(dx).
\end{eqnarray}
It is now seen that formally taking the limit as $T\to\infty$ we arrive at the right-hand side of (\ref{eq:3.71}). It remains to justify this procedure.

Denote 
$$U_1(T,s,y)=\int^{T^{1-s}-1}_0\int_{\errita^d}p_u(y-zT^{-s/\alpha})\varphi(z)
\chi((u+1)T^{s-1})dzdu$$
and
$$U_2(T,x)=\int^1_0\int_{\errita^d}p_1(xT^{-s/\alpha}-y)U^{1+\beta}_1(T,s,y)dyds.$$

We will need the following fact, which can be found, e.g., in \cite{14} 
(Lemma 5.3)
\vglue.3cm
\noindent
\begin{equation}
\label{eq:3.75}
\sup\limits_{x\in\errita^d}(1+|x|^{d-\alpha})|G\varphi(x)|<\infty,
\quad \varphi\in{\cal S}(\erre^d), d>\alpha .
\end{equation}
We have
\begin{eqnarray*}
U_1(T,s,y)&\leq&C\int_{\errita^d}\frac{1}{|y-zT^{-s/\alpha}|^{d-\alpha}}\varphi(z)dz\\
&=&\frac{C_1}{|y|^{d-\alpha}}|yT^{s/\alpha}|^{d-\alpha}G\varphi(yT^{s/\alpha})\\
&\leq&\frac{C_2}{|y|^{d-\alpha}},
\end{eqnarray*}
by (\ref{eq:3.75}). On the other hand, using the well known estimate
\begin{equation}
\label{eq:3.76}
p_1(x)\leq \frac{C}{1+|x|^{d+\alpha}},
\end{equation}
we have
$$\int^1_0p_1(xT^{-s/\alpha}-y)ds\leq C_3\frac{1+|x|^{d+\alpha}}{1+|y|^{d+\alpha}},
$$
hence it is not hard to see that $U_2(T,x)$
 converges pointwise as $T\to\infty$ and is bounded in $x,T$, since $(d-\alpha)(1+\beta)<d$ by (\ref{eq:2.19}).
 This proves (\ref{eq:3.71}) by (\ref{eq:3.73}) and (\ref{eq:3.74}).

Next observe that (\ref{eq:3.19}) implies that
\begin{equation}
\label{eq:3.77}
I_2(T)\leq C\frac{H_T}{F^2_T}\int_{\errita^d}G(\varphi G\varphi)(x)\mu(dx),
\end{equation}
hence (\ref{eq:3.17}) follows by (\ref{eq:3.72}) and (\ref{eq:2.20}).

Using (\ref{eq:3.72}), by (\ref{eq:3.21}) we have
\begin{eqnarray*}
J_1(T)&\leq&\frac{C}{F^{1+\beta}_T}\frac{H_T}{F^{1+\beta}_T}\int_{\errita^d}\int^T_0{\cal T}_{T-s}\biggl(\int^s_0{\cal T}_{s-u}\varphi du\biggr)^{1+\beta}(x)ds\mu(dx)\\
&\leq&\frac{C_1}{F^{1+\beta}_T}\to0,
\end{eqnarray*}
since from the proof of (\ref{eq:3.71}) it follows that
\begin{equation}
\label{eq:3.78}
\sup_{T\geq 2}\sup_{x\in\errita^d}\frac{1}{\log T}\int^T_0{\cal T}_{T-s}\biggl(\int^s_0{\cal T}_{s-u}\varphi du\biggr)^{1+\beta}(x)ds<\infty.
\end{equation}

To prove (\ref{eq:3.24}) it suffices to note that by (\ref{eq:3.22}) and 
(\ref{eq:3.78}),
$$J_2(T)\leq C\frac{H_T}{F^{(1+\beta)^2}_T}
T(\log T)^{1+\beta}=C\frac{T}{H^\beta_T}\to 0,$$
by (\ref{eq:2.16}).
 This completes the proof of (\ref{eq:3.5}).

To show tightness we prove (\ref{eq:3.33}) and (\ref{eq:3.34}) with $h=1$ and $\sigma$ satisfying (\ref{eq:3.45})
 (such $\sigma$ exists by (\ref{eq:2.19})). Recall that now we consider $t_1,t_2$  such that $0<\varepsilon <t_1<t_2\leq 1$, hence in (\ref{eq:3.61}) the integral on $s$ is taken over $[\varepsilon, 1]$. In  (\ref{eq:3.61}) we estimate $|\widehat{\varphi}(x)\widehat{\mu}(x+y)|$ by a constant and we integrate with respect to $x$, obtaining
$$A(T)\leq C\frac{H_TT^{2-d/\alpha}}{F^2_T}\int^1_\varepsilon\int^1_s\int_{\errita^d}s^{-d/\alpha}|\widehat{\varphi}(y)|e^{-T(u-s)|y|^\alpha}\chi(u)\chi(s)dyduds.$$
By (\ref{eq:3.47}) and(\ref{eq:2.20})
  we have
\begin{eqnarray*}
A(T)&\leq&C_1\varepsilon^{-d/\alpha}T^{1-d/\alpha+\sigma}(t_2-t_1)^\sigma\int^1_\varepsilon \chi(s)ds\int_{\errita^d}|\widehat{\varphi}(y)||y|^{\alpha(\sigma-1)}dy\\
&\leq&C_2(\varepsilon)T^{1-d/\alpha+\sigma}(t_2-t_1)^{1+\sigma}.
\end{eqnarray*}
Hence (\ref{eq:3.33}) follows by (\ref{eq:3.36}), (\ref{eq:3.45}) and (\ref{eq:2.19}).

Now we pass to the proof of (\ref{eq:3.34}). In this case the formula 
(\ref{eq:3.51}) has the form
\begin{equation}
\label{eq:3.79}
I_1(T)=Q_1(T)+Q_2(T),
\end{equation}
where
\begin{eqnarray}
\label{eq:3.80}
Q_1(T)&=&\frac{1}{\log T}\int_{\errita^d}\int^{\varepsilon/2}_0
p_s(xT^{-1/\alpha}-y)\biggl(\int^1_s\int_{\errita^d}p_{u-s}(y-z)\widetilde{\varphi}_T(z)\chi(u)dzdu\biggr)^{1+\beta}dyds\mu(dx),\phantom{dddddd}\\
\label{eq:3.81}
Q_2(T)&=&\frac{1}{\log T}\int_{\errita^d}\int^1_{\varepsilon/2}\int_{\errita^d}\ldots dyds\mu(dx).
\end{eqnarray}
In (\ref{eq:3.80}) we have $u-s>\varepsilon/2$, hence $p_{u-s}(y-z)\leq C(\varepsilon/2)^{-d/\alpha}$. Therefore
\begin{equation}
\label{eq:3.82}
Q_1(T)\leq C_1(\varepsilon)\biggl(\int_{\errita^d}\varphi(z)dz\biggr)^{1+\beta}\mu(\erre^d)(t_2-t_1)^{1+\beta}\leq C_2(\varepsilon)(t_2-t_1)^{1+\sigma}.
\end{equation}

In (\ref{eq:3.81}) we estimate $p_s(xT^{-1/\alpha}-y)$ by 
$C({\varepsilon}/{2})^{-d/\alpha}$, hence
$$Q_2(T)\leq C_3(\varepsilon)\int^1_0\int_{\errita^d}\biggl(\int^1_s\int_{\errita^d}p_{u-s}(y-z)\widetilde{\varphi}_T(z)dz \chi(u)du\biggr)^{1+\beta}dyds.$$
The last expression is identical with the estimate of $I\!\!I$ in \cite{6}, and it was shown there that it can be estimated by $C(t_2-t_1)^{1+\sigma}$, provided that $d<\alpha(1+\beta)/\beta$, which holds in our case by 
(\ref{eq:2.19}). This, together with (\ref{eq:3.82}), (\ref{eq:3.79}) and 
(\ref{eq:3.37}), proves (\ref{eq:3.34}), so the proof of part (b) is complete.

{\it Proof of part (c).}
First we show that under the assumption (\ref{eq:2.21}),
\begin{equation}
\label{eq:3.83}
\sup_{x\in\errita^d}G((G\varphi)^{1+\beta})(x)<\infty, \quad \varphi\in{\cal S}(\erre^d).
\end{equation}
We have
$$
G(G\varphi)^{1+\beta}(x)=C_{\alpha,d}\int_{|x-y|<1}\frac{1}{|x-y|^{d-\alpha}}(G\varphi)^{1+\beta}(y)dy
+C_{\alpha,d}g*(G\varphi)^{1+\beta}(x),
$$
where $g(x)=\UNO_{[1,\infty)}(|x|)|x|^{\alpha-d}$. The first term is bounded since $G\varphi$ is bounded. To show that the second term is bounded it suffices to find $p,q\geq 1, 1/p+1/q=1$, such that $g\in L^p$ and $(G\varphi)^{1+\beta}\in L^q$. Fix $q$ such that
$$\frac{d}{\alpha}>q>\max\left\{\frac{d}{(1+\beta)(d-\alpha)},1\right\}$$
(such $q$ exists by  (\ref{eq:2.21})). Then  (\ref{eq:3.75}) implies that 
$G^{1+\beta}\varphi\in L^q$, and it is clear that $g\in L^p$ for the corresponding $p$.

We now study the convergence of $I_1,I_2$ and $I_3$ defined by (\ref{eq:3.13})-(\ref{eq:3.15}). We have
\begin{equation}
\label{eq:3.84}
I_1(T)=\int_{\errita^d}\int^T_0{\cal T}_s\biggl(\int^{T-s}_0{\cal T}_u\varphi 
\chi\left(\frac{s}{T}+\frac{u}{T}\right)du\biggr)^{1+\beta}(x)ds\mu(dx).
\end{equation}
It is not difficult to see that (\ref{eq:3.83}) implies that
\begin{equation}
\label{eq:3.85}
\lim_{T\to\infty}I_1(T)=\int_{\errita^d}
G((G\varphi)^{1+\beta})(x)\mu(dx)(\chi(0))^{1+\beta}.
\end{equation}
By (\ref{eq:3.21}), (\ref{eq:3.72}) and (\ref{eq:3.83}) we have
$$J_1(T)\leq \frac{C}{H_T}\int_{\errita^d}G((G\varphi)^{1+\beta})(x)\mu(dx)\to 0.$$
Similarly, using (\ref{eq:3.22}) and (\ref{eq:3.83}),
$$J_2(T)\leq C\frac{T}{H_T^\beta}\to 0,$$
by (\ref{eq:2.16}).
 Hence we obtain (\ref{eq:3.18}), by (\ref{eq:3.20}) and 
(\ref{eq:3.85}). Next, for $\beta<1$,  (\ref{eq:3.77}) implies (\ref{eq:3.17}). This together with (\ref{eq:3.18}), (\ref{eq:3.85}) and (\ref{eq:3.12}) yield (\ref{eq:3.5}) in the case $\beta<1$ with $X$ determined by (\ref{eq:2.23}). To obtain 
(\ref{eq:3.5}) in the case $\beta=1$ it remains to show that
\begin{equation}
\label{eq:3.86}
\lim_{T\to\infty}I_2(T)=\int_{\errita^d}G(\varphi G\varphi)\mu(dx)
(\chi(0))^2.
\end{equation}
Using (\ref{eq:3.14}) and (\ref{eq:3.8}) we write
$$I_2(T)=I'_2(T)-I''_2(T)-I'''_3(T),$$
where
\begin{eqnarray*}
I'_2(T)&=&\int_{\errita^d}\int^T_0{\cal T}_{T-s}\left(\varphi \chi
\left(\frac{T-s}{T}\right)\int^s_0{\cal T}_{s-u}\varphi 
\chi\left(\frac{T-s}{T}\right)du\right)(x)ds\mu(dx),\\
I''_2(T)&=&\int_{\errita^d}\int^T_0{\cal T}_{T-s}\left(\varphi 
\chi\left(\frac{T-s}{T}\right)\int^s_0{\cal T}_{s-u}\varphi 
\chi\left(\frac{T-u}{T}\right)v_T(\cdot,u)du\right)(x)ds\mu
(dx),\\
I'''_2(T)&=&\frac{V}{2}H^{1/2}_T\int_{\errita^d}\int^T_0{\cal T}_{T-s}\left(\varphi \chi\left(\frac{T-s}{T}\right)\int^s_0{\cal T}_{s-u}v^2_T(\cdot, u)du\right)(x)ds\mu(dx).
\end{eqnarray*}
It is easy to see that $I'_2(T)$ converges to the right-hand side of 
(\ref{eq:3.86}). To show that $I''_2(T)$ and $I'''_2(T)$ converge to $0$, we first apply (\ref{eq:3.10}), and then use (\ref{eq:3.72}) and (\ref{eq:3.83}).

Finally, we pass to the proof of tightness. For $0<\varepsilon\leq t_1<t_2$, by (\ref{eq:3.38a}) and (\ref{eq:3.25})
 we have
\begin{eqnarray*}
A(T)&\leq&C\mu(\erre^d)\sup_{x\in\errita^d}\int^T_0{\cal T}_s\left(\varphi\int^T_s{\cal T}_{u-s}\varphi 
\chi\left(\frac{u}{T}\right)du\right)(x)\chi\left(\frac{s}{T}\right)ds\\
&\leq&C_1\sup_{x\in\errita^d}\int^{t_2T}_{t_1T}\int_{\errita^d}s^{-d/\alpha}p_1((x-y)x^{-1/\alpha})\varphi(y)
\int^{t_2T-s}_{0}{\cal T}_u\varphi(y)dudyds\\
&\leq&C_2\varepsilon^{-d/\alpha}T^{1-d/\alpha}(t_2-t_1)\int_{\errita^d}\varphi(y)dy\sup_y \int^{(t_2-t_1)T}_0{\cal T}_u\varphi(y)du\\
&\leq&C_3(\varepsilon)T^{1-d/\alpha+\sigma}(t_2-t_1)^{1+\sigma}(\sup_y G\varphi(y))^{1-\sigma}\\
&\leq&C_4(\varepsilon)(t_2-t_1)^{1+\sigma},
\end{eqnarray*}
for any
\begin{equation}
\label{eq:3.87}
0<\sigma<\left(\frac{d}{\alpha}-1\right)\wedge1,
\end{equation}
so we obtain (\ref{eq:3.33}).

To derive (\ref{eq:3.34}) we use (\ref{eq:3.37}),  (\ref{eq:3.84}) and 
(\ref{eq:3.25}),
obtaining
\begin{equation}
\label{eq:3.88}
I_1(T)\leq \mu(\erre^d)\sup_{x\in\errita^d}(Z_1(T,x)+Z_2(T,x)+Z_3(T,x)),
\end{equation}
where
\begin{eqnarray}
\label{eq:3.89}
Z_1(T,x)&=&\int^{t_1T/2}_0{\cal T}_s\biggl(\int^{t_2T-s}_{t_1T-s}{\cal T}_u\varphi du\biggr)^{1+\beta}(x)ds,\\
\label{eq:3.90}
Z_2(T,x)&=&\int^{t_1T}_{t_1T/2}{\cal T}_s\biggl(\int^{t_2T-s}_{t_1T-s}{\cal T}_u\varphi du\biggr)^{1+\beta}(x)ds,\\
\label{eq:3.91}
Z_3(T,x)&=&\int^{t_2T}_{t_1T}{\cal T}_s\biggl(\int^{t_2T-s}_0{\cal T}_u \varphi du\biggr)^{1+\beta}(x)ds.
\end{eqnarray}

By self-similarity we have
$$
Z_1(T,x)\leq C\int^{t_1T/2}_0\biggl(\int^{t_2T-s}_{t_1T-s}u^{-d/\alpha}du\biggr)^{1+\beta}ds.$$
As $u\geq t_1T-s\geq t_1T/2\geq \varepsilon T/2$, we get
\begin{eqnarray}
Z_1(T,x)&\leq&C_1\varepsilon^{1-(d/\alpha)(1+\beta)}T^{2+\beta-(d/\alpha)(1+\beta)}(t_2-t_1)^{1+\beta},\nonumber\\
\label{eq:3.92}
&\leq&C_2(\varepsilon)(t_2-t_1)^{1+\beta},
\end{eqnarray}
by  (\ref{eq:2.21}).

To estimate $Z_2$ we first use the bound 
$p_s(x-y)\leq C({T\varepsilon}/{2})^{-d/\alpha}$ for $s\geq t_1T/2$. After obvious substitutions we have
\begin{equation}
\label{eq:3.93}
Z_2(T,x)\leq C(\varepsilon)(Z'_2(T,x)+Z''_2(T,x)),
\end{equation}
where
\begin{eqnarray}
\label{eq:3.94}
Z'_2(T,x)&=&T^{-d/\alpha}\int^1_0\int_{\errita^d}\left(\int^{(t_2-t_1)T+s}_s{\cal T}_u\varphi(y)du\right)^{1+\beta}dyds,\\
\label{eq:3.95}
Z''_2(T,x)&=&T^{-d/\alpha}\biggl(\int^{t_1T/2}_1\int_{\errita^d}\ldots dyds\biggr)^+.
\end{eqnarray}
For any $0<\sigma\leq\beta$ we have
\begin{eqnarray}
\biggl(\int^{(t_2-t_1)T+s}_s{\cal T}_u\varphi(y)du\biggr)^\beta&\leq&(G\varphi(y))^{\beta-\sigma}\biggl(\sup_{y\in\errita^d}\varphi(y)\biggr)^\sigma((t_2-t_1)T)^\sigma.\nonumber\\
\label{eq:3.96}
&\leq&C(t_2-t_1)^\sigma T^\sigma,
\end{eqnarray}
by (\ref{eq:3.72}). Applying this to (\ref{eq:3.94})
we obtain 
\begin{equation}
\label{eq:3.97}
Z'_2(T,x)\leq C_1T^{-d/\alpha+1+\sigma}(t_2-t_1)^{1+\sigma}
\leq C_1(t_2-t_1)^{1+\sigma},
\end{equation}
provided that
\begin{equation}
\label{eq:3.98}
0<\sigma<\left(\frac{d}{\alpha}-1\right)\wedge \beta.
\end{equation}
In order to estimate $Z''_2$ we notice that for $d>\alpha$ and $0<a<b$,
$$\int^b_a{\cal T}_u\varphi(y)du\leq C\int^b_au^{-d/\alpha}du\leq
\left\{\begin{array}{l}
C(b-a)a^{-d/\alpha},\\
C_1a^{1-d/\alpha}.\end{array}\right.
$$
Using these two bounds, instead of (\ref{eq:3.96}) we now have for $0<\sigma\leq\beta$,
$$\biggl(\int^{(t_2-t_1)T+s}_s{\cal T}_u\varphi(y)du\biggr)^\beta\leq C_2s^{(1-d/\alpha)(\beta-\sigma)}((t_2-t_1)Ts^{-d/\alpha})^\sigma.$$
Putting this into (\ref{eq:3.95}) we obtain for $t_1T/2 >1$,
\begin{eqnarray}
Z''_2(T,x)&\leq&C_3T^{-d/\alpha+1+\sigma}\int^{t_1T/2}_1s^{-\sigma+(1-d/\alpha)\beta}dx (t_2-t_1)^{1+\sigma}\nonumber\\
&\leq&C_3T^{-d/\alpha+1+\sigma}\max(1,T^{1-\sigma+(1-d/\alpha)\beta}\log T)(t_2-t_1)^{1+\sigma}\nonumber\\
\label{eq:3.99}
&\leq&C_4(t_2-t_1)^{1+\sigma},
\end{eqnarray}
provided that (\ref{eq:3.98}) holds, and we also use (\ref{eq:2.21}) .
Combining (\ref{eq:3.93}), (\ref{eq:3.97}) and (\ref{eq:3.99}) we arrive at
\begin{equation}Z_2(T,x)\leq C(t_2-t_1)^{1+\sigma}\label{eq:3.100a}\end{equation}
for $\sigma$ satisfying (\ref{eq:3.98}).

Finally, by the H\"older inequality and using the fact that $t_2T-s\leq(t_2-t_1)T$, we have
\begin{eqnarray}
Z_3(T,x)&\leq&\biggl(\int^{t_2T}_{t_1T}\int_{\errita^d}p_s(x-y)dyds\biggr)^{1/(2+\beta)}\nonumber\\
&&\times\biggl[\int^{t_2T}_{t_1T}\int_{\errita^d}p_s(x-y)\left(\int^{(t_2-t_1)T}_0{\cal T}_u\varphi(y)du\right)^{2+\beta}dyds\biggr]^{(1+\beta)/(2+\beta)}\nonumber\\
&\leq&C((t_2-t_1)T)^{1/(2+\beta)}\biggl[\int^{t_2T}_{t_1T}\int_{\errita^d}s^{-d/\alpha}(G\varphi(y))^{2+\beta-\sigma}((t_2-t_1)T)^\sigma dyds\biggr]^{(1+\beta)/(2+\beta)}\nonumber\\
\label{eq:3.100}
&&
\end{eqnarray}
for any $0<\sigma\leq 2+\beta$, by an argument as in (\ref{eq:3.96}). Observe that by
 (\ref{eq:3.75}),
$$\int_{\errita^d}(G\varphi(y))^{2+\beta-\sigma}dy<\infty$$
for $\sigma$ sufficiently small, satisfying
\begin{equation}
\label{eq:3.101}
\frac{d}{\alpha}>\frac{2+\beta-\sigma}{1+\beta-\sigma}.
\end{equation}

Hence, by (\ref{eq:3.100}) we have
\begin{eqnarray}
Z_3(T,x)&\leq&C(\varepsilon)(t_2-t_1)^{1+\sigma(1+\beta)/(2+\beta)}T^{1+\sigma(1+\beta)/(2+\beta)-(d/\alpha)(1+\beta)/(2+\beta)}\nonumber\\
\label{eq:3.102}
&\leq&C(\varepsilon)(t_2-t_1)^{1+\sigma(1+\beta)/(2+\beta)},
\end{eqnarray}
provided that
\begin{equation}
\label{eq:3.103}
\sigma<\frac{d}{\alpha}\frac{1+\beta}{2+\beta}-1.
\end{equation}
Combining (\ref{eq:3.88}), (\ref{eq:3.92}), (\ref{eq:3.100a})
 and (\ref{eq:3.102}),
we conclude that (\ref{eq:3.34}) holds (with 
$\sigma(1+\beta)/(2+\beta)$ instead of $\sigma$) for any $\sigma$ 
satisfying (\ref{eq:3.98}), (\ref{eq:3.101}) and
(\ref{eq:3.103}).

The proof of Theorem 2.7 is complete. \hfill$\Box$
\vglue.5cm
\noindent
{\bf Proof of Proposition 2.9.} Only part (d) of the proposition needs to be proved. The argument is similar to that used in the proof of Theorem 2.7 in \cite{6}.

Observe that the finite-dimensional distributions of the process $\zeta$ defined by (\ref{eq:2.2}) are determined by
\begin{eqnarray}
\lefteqn{E{\rm exp}
\{
i(z_1\xi_{t_1}+\ldots+z_k\xi_{t_k})
\}
}\nonumber\\
&=&{\rm exp}
\biggl\{
-\int_{\errita^{d+1}}
\biggl[
\biggl|
\sum^k_{j=1}z_jp_r^{1/(1+\beta)}(x)
\UNO_{[0,t_j]}(r)\int^{t_j}_rp_{u-r}(x)du
\biggr|
^{1+\beta}
\biggr.
\nonumber\\
\label{eq:3.104}
&&\times\biggl.
\biggl(
1-i{\rm sgn}
\biggl(
\sum^k_{j=1}z_jp^{1/(1+\beta)}_r(x)\UNO_{[0,t_j]}(r)\int_r^{t_j}p_{u-r}(x)du
\biggr)
{\rm tan}\frac{\pi}{2}(1+\beta)
\biggr)
\biggr]
drdx
\end{eqnarray}
(see Proposition 3.4.2 of \cite{15}).

Denote
\begin{eqnarray*}
D^+_T&=&D_T(1,z;u,v,s,t),\quad z>0,\\
D^-_T&=&D_T(1,-z;u,v,s,t),\quad z>0,
\end{eqnarray*}
(see (\ref{eq:2.12})). It suffices to show that for fixed $0\leq u<v<s<t$ and 
$z>0$,
\begin{equation}
\label{eq:3.105}
D^+_T\leq CT^{-d/\alpha},\quad D^-_T\leq CT^{-d/\alpha},
\end{equation}
and for $T$ sufficiently large,
\begin{equation}
\label{eq:3.106}
D^+_T\geq CT^{-d/\alpha}
\end{equation}
(see (\ref{eq:2.11})).

It will be convenient to denote
\begin{eqnarray*}
f=f(x,r)&=&z\int^{t+T}_{s+T}p_{r'-r}(x)dr',\\
g_1=g_1(x,r)&=&\int^v_up_{r'-r}(x)dr',\\
g_2=g_2(x,r)&=&\int^v_rp_{r'-r}(x)dr'.
\end{eqnarray*}

It is not difficult to see that by (\ref{eq:3.104}),
\begin{eqnarray}
D^+_T&=&C\left[\int^u_0\int_{\errita^d}p_r(x)((f+g_1)^{1+\beta}-f^{1+\beta}-g^{1+\beta}_1)dxdr\right.\nonumber\\
\label{eq:3.107}
&&\left.\qquad\qquad +\int^v_u\int_{\errita^d}p_r(x)((f+g_2)^{1+\beta}-f^{1+\beta}-g^{1+\beta}_2)dxdr\right].
\end{eqnarray}
By the elementary inequality
$$0\leq(a+b)^{1+\beta}-a^{1+\beta}-b^{1+\beta}\leq (1+\beta)ab^\beta,\quad a,b\geq 0,\, 0<\beta\leq 1,$$
and  the estimate
$$f(x,r)\leq CT^{-d/\alpha},$$
we have
\begin{eqnarray}
\label{eq:3.108}
D^+_T&\leq&C_1\int_{\errita^d}\left[\int^u_0p_r(x)fg^\beta_1dr+\int^v_up_r(x)fg^\beta_2dr\right]dx\\
&\leq&C_2T^{-d/\alpha}\int_{\errita^d}\left(\int^v_0p_r(x)dr\right)^{1+\beta}dx\nonumber\\
&\leq&C_3T^{-d/\alpha},\nonumber
\end{eqnarray}
by (\ref{eq:3.3}). One can show that for $D^-_T$ the estimate (\ref{eq:3.108})
also holds (see \cite{6} for details).
Hence (\ref{eq:3.105}) follows.

Next, by (\ref{eq:3.107}),
\begin{eqnarray*}
D^+_T&\geq&C\int^{(u+v)/2}_u\int_{|x|\leq 1}p_r(x)((f+g_2)^{1+\beta}-f^{1+\beta}-g^{1+\beta})dxdr\\
&\geq&C_1\int^{(u+v)/2}_u\int_{|x|\leq 1}((f+g_2)^{1+\beta}-f^{1+\beta}-g^{1+\beta})dxdv,
\end{eqnarray*}
and this is exactly the right-hand side of (4.18) in \cite{6}, and it was proved there that it is greater than $CT^{-d/\alpha}$ for large $T$.
Thus (\ref{eq:3.106}) holds. \hfill$\Box$
\vglue.5cm
\noindent
{\bf Proof of Theorem 2.10.} The theorem can be proved using the corresponding version of the general scheme (see (\ref{eq:3.38})
 and the discussion following it). The arguments are similar to those carried out in the branching case and they are easier, therefore we omit the proof. We only indicate how to obtain the process $\zeta$ in  part (a).

It is easy to see that $I\!\!I_1(T)$ defined by (\ref{eq:3.39}) can be written as
\begin{eqnarray*}
 I\!\!I_1(T)&=&\int_{\errita^d}\int^1_0\int_{\errita^d}\int^1_s\int_{\errita^d}s^{-d/\alpha}p_1((x-y)T^{-1/\alpha}s^{-1/\alpha})\varphi(y)
\chi(s)(u-s)^{-d/\alpha}\\
&&\times p_1((y-z)T^{-1/\alpha}(u-s)^{-1/\alpha})
\varphi(z)\chi(u)dzdudyds\mu(dx)\\
&\to&p^2_1(0)\mu(\erre^d)\left(\int_{\errita^d}\varphi(y)dy\right)^2\int^1_0\int^1_s(u-s)^{-d/\alpha}s^{-d/\alpha}\chi(s)\chi(u)duds,
\end{eqnarray*}
and this is exactly the logarithm of right-hand side of (\ref{eq:3.41}) with
$$X=\left(\frac{2p^2_1(0)\mu(\erre^d)}{1-\frac{d}{\alpha}}\right)^{1/2}
\lambda \rho,$$
and $\rho$ is  Gaussian with covariance (\ref{eq:2.26}). \hfill$\Box$
\vglue.5cm
\noindent
{\bf Acknowledgment.} We thank the Institute of Mathematics,  National University of Mexico (UNAM), where part of this work was done.

\end{document}